\documentclass[12pt,leqno]{article}
\usepackage{amsfonts}
\usepackage{amsmath, amsthm, amsfonts, amssymb, color}
\usepackage{mathrsfs}
\renewcommand{\bar}{\overline}
\renewcommand{\hat}{\widehat}
\renewcommand{\tilde}{\widetilde}

\theoremstyle{definition}

\definecolor{wco}{rgb}{0.5,0.2,0.3}

\numberwithin{equation}{section} \theoremstyle{remark}

\newcommand{\ua}{\uparrow}
\newcommand{\da}{\downarrow}

\title{{\bf Asymptotic behaviors for distribution dependent SDEs driven by fractional Brownian motions}
}
\author{
{\bf  Xiliang Fan$^{a)}$, Ting Yu$^{a)}$,  Chenggui Yuan$^{b)}$
 }\\
\footnotesize{$^{a)}$School of Mathematics and Statistics, Anhui Normal University, Wuhu 241002, China}\\
\footnotesize{$^{b)}$Department of Mathematics, Swansea University, Bay Campus, SA1 8EN, UK}\\
\footnotesize{\sf{fanxiliang0515@163.com},\ \sf{tingyu1328@163.com},\ \sf{c.yuan@swansea.ac.uk }}\\
}

\begin{document}
\def\R{\mathbb R}
\def\N{\mathbb N}
\def\E{\mathbb E}
 \def\H{\mathbb H}
\def\Q{\mathbb Q}
\def\P{\mathbb{P}}
\def\S{\mathbb{S}}
\def\Y{\mathbb{Y}}
\def\W{\mathbb{W}}
\def\D{\mathbb{D}}
\def\A{\mathcal{A}}
\def\F{\mathcal{F}}
\def\G{\mathcal{G}}
\def\H{\mathcal{H}}
\def\cS{\mathcal{S}}

\def\sA{\mathscr{A}}
 \def\sB{\mathscr {B}}
 \def\sC{\mathscr {C}}
  \def\sD{\mathscr {D}}
 \def\sF{\mathscr{F}}
\def\sG{\mathscr{G}}
\def\sL{\mathscr{L}}
\def\sP{\mathscr{P}}
\def\sM{\mathscr{M}}
\def\eq{\equation}
\def\beg{\begin}
\def\ep{\epsilon}
\def\ve{\varepsilon}
\def\vp{\varphi}
\def\vr{\varrho}
\def\om{\omega}
\def\Om{\Omega}
\def\si{\sigma}
\def\ff{\frac}
\def\sq{\sqrt}
\def\kk{\kappa}
\def\de{\delta}
\def\<{\langle} \def\>{\rangle}
\def\Ga{\Gamma}
\def\ga{\gamma}
\def\na{\nabla}
\def\be{\beta}
\def\al{\alpha}
\def\pp{\partial}
 \def\ti{\tilde}
\def\1{\lesssim}
\def\ra{\rightarrow}
\def\da{\downarrow}
\def\upa{\uparrow}
\def\l{\ell}
\def\8{\infty}
\def\3{\triangle}
 \def\DD{\Delta}
\def\m{{\bf m}}
\def\B{\mathbf B}
\def\e{\text{\rm{e}}}
\def\la{\lambda}
\def\th{\theta}

\def\d{\text{\rm{d}}}
   \def\ess{\text{\rm{ess}}}
\def\Ric{\text{\rm{Ric}}} \def \Hess{\text{\rm{Hess}}}
 \def\ua{\underline a}
  \def\Ric{\text{\rm{Ric}}}
\def\cut{\text{\rm{cut}}}      \def\alphaa{\mathbf{r}}     \def\r{r}
\def\gap{\text{\rm{gap}}} \def\prr{\pi_{{\bf m},\varrho}}  \def\r{\mathbf r}

\def\Tilde{\tilde} \def\TILDE{\tilde}\def\II
{\mathbb I}
\def\i{{\rm in}}\def\Sect{{\rm Sect}}

\renewcommand{\bar}{\overline}
\renewcommand{\hat}{\widehat}
\renewcommand{\tilde}{\widetilde}

\allowdisplaybreaks
\maketitle
\begin{abstract}
In this paper, we study small-time asymptotic behaviors for a class of distribution dependent stochastic differential equations
driven by fractional Brownian motions with Hurst parameter $H\in(1/2,1)$ and magnitude $\ep^H$.
By building up a variational framework and two weak convergence criteria in the factional Brownian motion setting,
we establish the large and moderate deviation principles for this type equations.
Besides, we also obtain the central limit theorem, in which the limit process solves a linear equation involving the Lions derivative of
the drift coefficient.
\end{abstract}
AMS Subject Classification: 60H10, 60G22

Keywords: Distribution dependent SDE; fractional Brownian motion; large deviation principle; moderate deviation principle; central limit theorem

\section{Introduction}

In this article, we study the asymptotic behaviors (including the large and moderate deviation principles and the central limit theorem) for distribution dependent stochastic differential equations (DDSDEs) in $\R^d$ with small fractional noises as follows:
\begin{align}\label{1-Int}
\d X_t^\epsilon =b(t,X_t^\ep,\sL_{X_t^\ep})\d t+\ep^H\si(t,\sL_{X_t^\ep})\d B_t^H,\ \ X_0^\ep =x,
\end{align}
where $\sL_{X_t^\ep}$ denotes the law of $X_t^\ep$, $\ep>0$ is a small parameter, $B^H$ is a fractional Brownian motion with Hurst parameter $H\in(1/2,1)$, the coefficients $b$ and $\si$ fulfill some appropriate conditions given in later sections.
Moreover, the integral with respect to $B^H$ is interpreted in the Wiener sense due to the determinacy of $\si(\cdot,\sL_{X_\cdot^\ep})$.
Before stating more precisely our approach and our results, let us recall a few results concerning distribution dependent SDEs and large and moderate deviation principles.

DDSDEs, also called McKean-Vlasov or mean-field SDEs, were first studied by McKean \cite{McKean66} to model plasma dynamics.
These DDSDEs describe limiting behaviours of individual particles in a large system of particles which interact in a mean-field sense, as the number of particles tends to infinity.
Another important feature of DDSDEs is their intrinsic link with nonlinear Fokker-Planck-Kolmogorov equation which may characterize the evolution of the marginal laws of DDSDEs.
For these reasons, DDSDEs have found applications in numerous fields such as statistical physics, mean-field games, mathematical finance and biology
(see, e.g., \cite{BT97,CD13,JW17,LL07} and references therein).
Recently, there has been a flourishing growth in the literature on DDSDEs,
one can refer to \cite{Buckdahn&L&Peng&ainer17a,Crisan&McMurray18a} for value functions and related PDEs,
\cite{HW,RW,Song,Wang18} for Harnack type inequalities, gradient estimates, Lions type derivative formulas, and many other aspects.

On the other hand, large and moderate deviation principles are to calculate the probability of a rare event, which investigate the asymptotic property of remote tails of a family of probability distribution.
In the case of stochastic processes, the idea lies in identifying a deterministic path around which the diffusion is concentrated with high probability,
which leads to a interpretation of the stochastic motion as a small perturbation of this deterministic path.
Since the original work of \cite{FW84} adopting method of contraction principle and the technique of time discretization,
large and moderate deviation principles for stochastic equations perturbed by Brownian motion or Poisson random measure have been extensively studied in the past several decades.
Another general approach for studying large and moderate deviation problems is the well-known weak convergence method introduced in \cite{BD00,BDM08},
which is based on a variational representation for positive functionals of Brownian motion or Poisson random measure.
After that, this approach has been wildly applied in various stochastic dynamical systems,
see, for example, \cite{BGJ17,BDG16,BDM11,DST19,DWZZ20,DXZZ17,HLL21,JY21,MSZ21,SY21,WZZ15,ZZ15} and the references therein.
In the study of moderate deviation principle (MDP), one is concerned with deviation probability of a lower order than that in large deviation principle (LDP), which actually bridges the gap between LDP and central limit theorem (CLT) (see more details in the next paragraph).

Contrary to the previously mentioned works, the aim of this paper is to study the asymptotic behaviors of equation \eqref{1-Int}
perturbed by a fractional motion $B^H$  with $H\in(1/2,1)$ as $\ep\ra0$.
Our primary motivation for doing this comes from the fact that the vast and growing literature on large and moderate deviation results for stochastic equations
focuses mostly on the noises of Markov property or independently identically distribution,
and results outside this class are few.
The only other works, to the best of our knowledge, are \cite{BS20} and \cite{BDS22} which treated classical SDEs driven by fractional Brownian motion (that is, distribution-free) with magnitude $\sqrt{\ep}$, rather than $\ep^H$, and investigated respectively LDP and MDP.
One possible reason for this is fractional Brownian motion $B^H$ with parameter $H\neq1/2$ is neither a Markov process nor a semimartingale
and then techniques based on the It\^{o} calculus are not applicable, substantial new difficulties will appear in this setting.
Recently, in \cite{FHSY} we proved the well-posedness of DDSDEs driven by fractional Brownian motions and then established the Bismut formulas for both non-degenerate and degenerate cases.
With the solution of \eqref{1-Int} in hand, a natural problem one can think of is the following:
can we obtain the small-time asymptotic results for this type equation?
The main purpose of this paper is to study the LDP, MDP and CLT of \eqref{1-Int} when $\ep\ra0$.
More precisely, let $X^0$ be the limit of $X^\ep$ in some sense, we are going to investigate the asymptotic behaviors for the path of the form
\begin{align*}
Y_t^\ep:=\ff {X_t^\ep-X_t^0}{\ep^H\kappa(\ep)},\ \ t\in[0,T].
\end{align*}

$\bullet$ In the case of the LDP, namely $\kappa(\ep)=1/\ep^H$, we show that $X^\ep$ satisfies the LDP with speed $\ep^{2H}$ (see Theorem \ref{Th(ldp)}).

$\bullet$ In the case of the CLT, namely $\kappa(\ep)=1$, we prove that as $\ep\ra0$, $\ff {X^\ep-X^0}{\ep^H}$ converges to a stochastic process which
solves a linear equation involving the Lions derivative of the coefficient $b$ (see Theorem \ref{Th(clt)}).

$\bullet$ In the case of the MDP, namely $\kappa(\ep)\ra\infty$ and $\ep^H\kappa(\ep)\ra0$ as $\ep\ra0$,
we derive that $Y^\ep$ satisfies LDP with speed $\kappa^{-2}(\ep)$ (see Theorem \ref{Th(mdp)}).\\
Here, let us point out that the MDP for $X^\ep$ refers to the LDP for $Y^\ep$ since the scaling by $\ep^H\kappa(\ep)$
means that the MDP is in the regime between the LDP and the CLT.

In order to obtain the large and moderate deviation principles, we shall adopt the weak convergence method.
With the help of a variational representation for random functional of fractional Brownian motion (see Lemma \ref{VaRe}),
we provide two sufficient conditions for the Laplace principle (which is equivalent to the LDP) to hold for functionals of fractional Brownian motion
(see Propositions \ref{Suf1(LDP)}, \ref{Suf2(LDP)}), which are extensions (or fractional versions) of \cite[Theorem 4.2]{BDM11}, \cite[Theorem 4.3]{LSZZ22} and \cite[Theorem 3.2]{MSZ21}.
Then we are allowed to verify the weak convergence properties of the controlled equations for DDSDEs concerned.
Let us stress that the identification of the correct controlled equations is a crucial step in the DDSDEs situation.

An interesting problem coming from the work is whether our results can be further improved in the sense of allowing the diffusion coefficient $\si$ to be dependent on $X_t$.
Our recent work \cite[Remark 3.2(ii)]{FHSY} suggests that there is an essential difficulty in establishing the well-posedness of multidimensional DDSDEs driven by fractional Brownian motions.
Recently, in the Brownian motion setting, some scholars have investigated the well-posendess of DDSDEs via the Girsanov transforms and the coupling argument (see, e.g., \cite{HW21b,W22}), which are available for fractional Brownian motion as well.
We hope that these tools are of help in future studies.
Another challenging question is whether there are similar results for DDSDEs driven by rough fractional Brownian motion, which corresponds to $H\in(0,1/2)$.
In the rough situation, the skeleton equation \eqref{SkEq} below might be ill-posedness, thus our techniques are currently not enough to give an affirmative answer.
We will leave this topic for future work.

The rest of the paper is organized as follows.
Section 2 is devoted to recalling some useful facts on fractional Brownian motion and the Lions derivative.
In Section 3, we present the framework of DDSDE driven by fractional Brownian motion and state some basic theory of the LDP.
Then we formulate the main results concerning the LDP, the MDP and the CLT for DDSDE driven by fractional Brownian motion.
In Section 4, we focus on proving our main results.
Section 5 will be devoted to the proofs of some auxiliary results.

\section{Preliminaries}

\subsection{Fractional Brownian motion}

In this part, we shall recall some important definitions and results concerning the fractional Brownian motion.
For a deeper discussion, we refer the reader to \cite{Alos&Mazet&Nualart01a,Biagini&Hu08a,Decreusefond&Ustunel98a,ND,Nualart&Saussereau09} and references therein.

For some fixed  $H\in(1/2,1)$.
we consider $(\Omega,\sF,\P)$ the canonical probability space associated with fractional Brownian motion with Hurst parameter $H$.
That is, $\Omega$ is the Banach space $C_0([0,T],\R^d)$ of continuous functions vanishing at $0$ equipped with the supremum norm,
$\sF$ is the Borel $\si$-algebra and $\P$ is the unique probability measure on $\Omega$ such that the canonical process $\{B^H_t; t\in[0,T]\}$ is a $d$-dimensional fractional Brownian motion with Hurst parameter $H$.
Recall that $B^H=(B^{H,1},\cdots,B^{H,d})$  is a centered Gaussian process, whose covariance structure is defined by
\begin{align*}
\E\left(B^{H,i}_tB^{H,j}_s\right)=R_H(t,s)\delta_{i,j},\ \ s,t\in[0,T],\  i,j=1,\cdots,d
\end{align*}
with  $R_H(t,s)=\frac{1}{2}(t^{2H}+s^{2H}-|t-s|^{2H})$.
Let $\{\sF_t\}_{t\in[0,T]}$ be the filtration generated by $B^H$.

We denote by $\mathscr{E}$ the set of step functions on $[0,T]$ with values in $\R^d$.
Let $\mathcal {H}$ be the Hilbert space defined as the closure of $\mathscr{E}$ with respect to the scalar product
\begin{align*}
\left\langle (\mathbb I_{[0,t_1]},\cdot\cdot\cdot,\mathbb I_{[0,t_d]}),(\mathbb I_{[0,s_1]},\cdot\cdot\cdot,\mathbb I_{[0,s_d]})\right\rangle_\H=\sum\limits_{i=1}^dR_H(t_i,s_i).
\end{align*}
Note that by \cite{Decreusefond&Ustunel98a,ND},  $R_H(t,s)$ has the following integral representation
\begin{align*}
 R_H(t,s)=\int_0^{t\wedge s}K_H(t,r)K_H(s,r)\d r,
\end{align*}
where $K_H(t,s)$ is the square integrable kernel given by
\begin{align*}
K_H(t,s)=C_Hs^{\ff 1 2-H}\int_s^t(r-s)^{H-\ff 3 2}r^{H-\ff 1 2}\d r, \ \ t>s
\end{align*}
with $C_H=\sqrt{\ff {H(2H-1)}{\mathcal{B}(2-2H,H-1/2)}}$ and $\mathcal{B}$ standing for the Beta function.
If $t\leq s$, we set $K_H(t,s)=0$.
The mapping $(\mathbb I_{[0,t_1]},\cdot\cdot\cdot,\mathbb I_{[0,t_d]})\mapsto\sum_{i=1}^dB_{t_i}^{H,i}$ can be extended to an isometry between $\H$ (also called the reproducing kernel Hilbert space) and the Gaussian space $\mathcal {H}_1$ associated to $B^H$. We denote this isometry by $\phi\mapsto B^H(\phi)$.

Now, let $(e_1,\cdots,e_d)$ designate the canonical basis of $\R^d$, one can introduce the linear operator $K_H^*:\mathscr{E}\rightarrow L^2([0,T],\R^d)$ defined by
\begin{align*}
K_H^*(\mathrm{I}_{[0,t]}e_i)=K_H(t,\cdot)e_i.
\end{align*}
By \cite{Alos&Mazet&Nualart01a}, it is known that the relation $\langle K_H^*\psi,K_H^*\phi\rangle_{L^2([0,T],\R^d)}=\langle\psi,\phi\rangle_\H$ holds for all $\psi,\phi\in\mathscr{E}$,
and then by the bounded linear transform theorem, $K_H^*$ can be extended to an isometry between $\mathcal{H}$ and $L^2([0,T],\R^d)$.
Consequently, by \cite{Alos&Mazet&Nualart01a} again,
there exists a $d$-dimensional Wiener process $W$ defined on $(\Omega,\sF,\P)$ such that $B^H$ has the following Volterra-type representation
\beg{align}\label{IRFor}
B_t^H=\int_0^tK_H(t,s)\d W_s, \ \ t\in[0,T],
\end{align}
and $\sF_t=\si\{\W_s:0\leq s\leq t\}$.
Moreover, let us point out that $K_H^*$ has the following representations: for any $\psi,\phi\in\H$,
\begin{align*}
(K_H^*\psi)(t)=\int_t^T\psi(s)\ff {\partial K_H(s,t)}{\partial s}\d s
\end{align*}
and
\begin{align}\label{Isom}
\langle K_H^*\psi,K_H^*\phi\rangle_{L^2([0,T],\R^d)}=\langle\psi,\phi\rangle_\H=H(2H-1)\int_0^T\int_0^T|t-s|^{2H-2}\langle\psi(s),\phi(t)\rangle_{\R^d}\d s\d t.
\end{align}
As a consequence, for any $\psi\in L^2([0,T],\R^d)$, one has
\begin{align}\label{EsH}
\|\psi\|^2_\H\leq2HT^{2H-1}\|\psi\|^2_{L^2}.
\end{align}
Besides, it can be shown that $L^{1/H}([0,T],\R^d)\subset\H$.

Next, we define the operator $K_H: L^2([0,T],\mathbb{R}^d)\rightarrow I_{0+}^{H+1/2}(L^2([0,T],\mathbb{R}^d))$ by
\begin{align*}
 (K_H f)(t)=\int_0^tK_H(t,s)f(s)\d s,
\end{align*}
where $ I_{0+}^{H+1/2 }$ is the left-sided  fractional Riemann-Liouville integral operator shown by \eqref{1-HLI} below with $\alpha=H+1/2$.
Let us mention that the space $I_{0+}^{H+1/2}(L^2([0,T],\mathbb{R}^d))$ is the fractional version of the Cameron-Martin space.

We remark that in the case of $H=1/2$ (usual Brownian motion), $K_H^*$ is the identity map on $L^2([0,T],\mathbb{R}^d)$,
$K_H(t,s)=\mathrm{I}_{[0,t]}(s)$ and $I_{0+}^{H+1/2}(L^2([0,T],\mathbb{R}^d))$ is the space of absolutely continuous functions, vanishing at zero, with a square integrable derivative.

Finally, we denote by $R_H=K_H\circ K_H^*:\mathcal{H}\rightarrow I_{0+}^{H+1/2}(L^2([0,T],\mathbb{R}^d))$ the operator
\begin{align}\label{RHop}
(R_H\psi)(t)=\int_0^tK_H(t,s)(K_H^*\psi)(s)\d s.
\end{align}
Since  $I_{0+}^{H+1/2}(L^2([0,T],\mathbb{R}^d))\subset C^{H}([0,T],\R^d)$ due to \cite[Theorem 3.6]{SKM93},
we know that for any $\psi\in\H$, $R_H\psi$ is H\"{o}lder continuous of order $H$, i.e.
\begin{align}\label{HORH}
R_H\psi\in C^{H}([0,T],\R^d),\ \ \psi\in\H.
\end{align}
On the other hand, using the Fubini theorem and the fact that $\ff {\partial K_H(s,r)}{\partial s}=C_H(\ff s r)^{H-\ff 1 2}(s-r)^{H-\ff 3 2}$,
it is easy to see that for every $\psi\in\H, R_H\psi$ is absolutely continuous and
\begin{align}\label{RFRH}
(R_H\psi)(t)=\int_0^t\left(\int_0^s\ff {\partial K_H}{\partial s}(s,r)(K_H^*\psi)(r)\d r\right)\d s.
\end{align}
We remark that the injection $R_H=K_H\circ K_H^*:\mathcal{H}\rightarrow\Omega$ embeds $\mathcal{H}$ densely into $\Omega$ and
for each $\psi\in\Omega^*\subset\H$ there holds
\begin{align*}
\E\e^{i\<B^H,\psi\>}=\exp\left(-\ff 1 2\|\psi\|^2_\H\right).
\end{align*}
Consequently, $(\Omega,\mathcal{H},\P)$ is an abstract Wiener space in the sense of Gross.
It is worth stressing that compared with the work \cite{Decreusefond&Ustunel98a} in which the authors have made another (but equivalent) choice for the underlying Hilbert space,
the choices of the Hilbert space and its embedding into $\Omega$ are not unique.

\subsection{The Lions derivative}

For later use, we introduce some basic facts about the Lions derivative.

For any $\th\in[1,\infty)$, $\sP_\th(\R^d)$ stands for the set of $\th$-integrable probability measures on $\R^d$,
and define the $L^\th$-Wasserstein distance on $\sP_\th(\R^d)$ as follows
\begin{align*}
\mathbb{W}_\th(\mu,\nu):=\inf_{\pi\in\sC(\mu,\nu)}\left(\int_{\R^d\times\R^d}|x-y|^\th\pi(\d x, \d y)\right)^\ff 1 \th,\ \ \mu,\nu\in\sP_\th(\R^d).
\end{align*}
Here $\sC(\mu,\nu)$ denotes the set of all probability measures on $\R^d\times\R^d$ with marginals $\mu$ and $\nu$.
It is well known that $(\sP_\th(\R^d),\mathbb{W}_\th)$ is a Polish space, usually referred to as the $\th$-Wasserstein space on $\R^d$.
Throughout this paper, let $|\cdot|$ and $\<\cdot,\cdot\>$  be the Euclidean norm and inner product, respectively,
and for a matrix, $\|\cdot\|$ denotes the operator norm.
The Lebesgue spaces  $L^2([0,T],\R^d)$ and  $ L^2(\R^d\ra\R^d,\mu)$ have the norms $\|\cdot\|_{L^2}$ and $\|\cdot\|_{L^2_\mu}$, respectively.
Let further  $\sL_X$ be the distribution of  random variable $X$.

\beg{defn}
Let $f:\sP_2(\R^d)\ra\R$ and $g:\R^d\times\sP_2(\R^d)\ra\R$.
\begin{enumerate}
\item[(1)]  $f$ is called $L$-differentiable at $\mu\in\sP_2(\R^d)$, if the functional
\begin{align*}
L^2(\R^d\ra\R^d,\mu)\ni\phi\mapsto f(\mu\circ(\mathrm{Id}+\phi)^{-1}))
\end{align*}
is Fr\'{e}chet differentiable at $0\in L^2(\R^d\ra\R^d,\mu)$. That is, there exists a unique $\gamma\in L^2(\R^d\ra\R^d,\mu)$ such that
\begin{align*}
\lim_{\|\phi\|_{L^2_\mu}\ra0}\ff{f(\mu\circ(\mathrm{Id}+\phi)^{-1})-f(\mu)-\mu(\<\gamma,\phi\>)}{\|\phi\|_{L^2_\mu}}=0.
\end{align*}
In this case, $\gamma$ is called the $L$-derivative of $f$ at $\mu$ and denoted by $D^Lf(\mu)$.

\item[(2)] $f$ is called $L$-differentiable on $\sP_2(\R^d)$, if the $L$-derivative $D^Lf(\mu)$ exists for all $\mu\in\sP_2(\R^d)$.
Furthermore, if for every $\mu\in\sP_2(\R^d)$ there exists a $\mu$-version $D^Lf(\mu)(\cdot)$ such that $D^Lf(\mu)(x)$ is jointly continuous in $(\mu,x)\in\sP_2(\R^d)\times\R^d$, we denote $f\in C^{(1,0)}(\sP_2(\R^d))$.

\item[(3)] $g$ is called differentiable on $\R^d\times\sP_2(\R^d)$, if for any $(x,\mu)\in\R^d\times\sP_2(\R^d)$,
$g(\cdot,\mu)$ is differentiable and $g(x, \cdot)$ is $L$-differentiable.
Furthermore, if $\nabla g(\cdot,\mu)(x)$ and $D^Lg(x,\cdot)(\mu)(y)$ are jointly continuous in $(x,y,\mu)\in\R^d\times\R^d\times\sP_2(\R^d)$,
we denote $g\in C^{1,(1,0)}(\R^d\times\sP_2(\R^d))$.

\end{enumerate}
\end{defn}

For a vector-valued function $f=(f_i)$ or a matrix-valued function $f=(f_{ij})$ with $L$-differentiable components, we simply write
\begin{align*}
D^Lf(\mu)=(D^Lf_i(\mu))  \ \  \mathrm{or}\ \ D^Lf(\mu)=(D^Lf_{ij}(\mu)).
\end{align*}

Besides, by  \cite[Theorem 6.5]{Cardaliaguet13} and \cite[Proposition 3.1]{RW}, we have the following useful formula for the $L$-derivative.
\beg{lem}\label{FoLD}
Let $(\Omega,\sF,\P)$ be an atomless probability space and $X,Y\in L^2(\Omega\ra\R^d,\P)$.
If $f\in C^{1,0}(\sP_2(\R^d))$, then
\begin{align*}
\lim_{\ve\da0}\ff {f(\sL_{X+\ve Y})-f(\sL_X)} \ve=\E\<D^Lf(\sL_X)(X),Y\>.
\end{align*}
\end{lem}

\subsection{Notations}

$C_b(\mathscr{E})$ denotes the set of all bounded continuous functions $f:\mathscr{E}\ra\R$ with the norm $\|f\|_\infty:=\sup_{x\in\mathscr{E}}|f(x)|$,
where $\mathscr{E}$ is a Polish space with the Borel $\sigma$-field $\mathcal{B}(\mathscr{E})$.
Let
\begin{align*}
\A=\left\{\phi:\phi\  \mathrm{is}\ \R^d\textit{-}\mathrm{valued}\ \F_t\textit{-}\mathrm{predictable\ process\ and \ \|\phi\|_\H^2<\infty\ \P\textit{-}\mathrm{a.s.} } \right\},
\end{align*}
and for each $M>0$, let
\begin{align*}
S_M=\left\{h\in\H: \ff 1 2\|h\|_\H^2\leq M\right\}.
\end{align*}
It is obvious that $S_M$ endowed with the weak topology is a Polish space.
Besides, define
\begin{align*}
\mathcal{A}_M:=\{\phi\in\mathcal{A}: \phi(\omega)\in S_M, \ \mathbb{P}\textit{-}\mathrm{a.s.}\}.
\end{align*}

\section{Framework and main results}

The main objective of this section concerns the problem of asymptotic behaviors for DDSDEs driven by fractional Brownian motions.
We first introduce the framework of DDSDEs driven by fractional Brownian motions.
Then we recall the definitions of large deviation principle and Laplace principle and present their relations.
Moreover, we provide the weak criteria for DDSDEs driven by fractional Brownian motions and formulate the main results including
large deviation principle, moderate deviation principle and central limit theorem for this type of equations.

\subsection{Framework}

We now consider the following distribution dependent SDE:
\begin{align}\label{DDsde}
\d X_t=b(t,X_t,\mathscr{L}_{X_t})\d t+\sigma(t,\mathscr{L}_{X_t})\d B_t^H,\ \ X_0=x\in \R^d , \ \ t\in[0,T].
\end{align}
Here $b:[0,T]\times \R^d\times\sP_\th(\R^d)\rightarrow\R^d$ and $\sigma:[0,T]\times\sP_\th(\R^d)\rightarrow\R^d\otimes\R^d$
with $\th\in[1,2]$ are measurable mappings, $B^H$ is a $d$-dimensional fractional Brownian motion with $H\in(1/2,1)$.
We assume that $b$ and $\si$ satisfy the following conditions:
\begin{enumerate}

\item[\textsc{\textbf{(H1)}}] There is a non-decreasing function $K(t)$ such that for any $t\in[0,T], x,y\in\R^d, \mu,\nu\in \sP_\th(\R^d)$,
\begin{align*}
|b(t,x,\mu)-b(t,y,\nu)|\leq K(t)(|x-y|+\mathbb{W}_\th(\mu,\nu)), \ \
\|\sigma(t,\mu)-\sigma(t,\nu)\|\leq K(t)\mathbb{W}_\th(\mu,\nu),
\end{align*}
and
\begin{align*}
|b(t,0,\delta_0)+\|\sigma(t,\delta_0)\|\leq K(t).
\end{align*}

\end{enumerate}
For each $p\geq1$, $\cS^p([0,T])$ denotes the space of  $\R^d$-valued, continuous $(\sF_t)_{t\in[0,T]}$-adapted processes $\phi$ on $[0,T]$ satisfying
\begin{align*}
\|\phi\|_{\cS^p}:=\bigg(\E\sup_{t\in[0,T]}|\phi_t|^p\bigg)^{1/p}<\infty,
\end{align*}
and $C_{p_1,p_2,\cdots}$ denotes generic constants, whose values may change from line to line and depend only on $p_1,p_2,\cdots$.

\beg{defn}\label{DeEq}
A stochastic process $X=(X_t)_{0\leq t \leq T}$ on $\R^d$ is called a solution of \eqref{DDsde},
if $X\in\cS^p([0,T])$ and $\P$-a.s.,
\begin{align*}
X_t=x+\int_0^tb(s,X_s,\sL_{X_s})\d s+\int_0^t\si(s,\sL_{X_s})\d B_s^H,\ \ t\in[0,T].
\end{align*}
\end{defn}
We remark here that $\int_0^t\si(s,\sL_{X_s})\d B_s^H$ can be interpreted as Wiener integral with respect to fractional Brownian motion
since $\si(\cdot,\sL_{X_\cdot})$ is a deterministic function.
According to \cite[Theorem 3.1]{FHSY},  \textsc{\textbf{(H1)}} implies that equation \eqref{DDsde} admits a unique solution $X\in\cS^p([0,T])$ with any $p\geq\th$ and $p>1/H$.

In order to investigate the asymptotic behaviors for the equation \eqref{DDsde},
for any fixed $\mu.\in C([0,T]; \sP_2(\R^d))$, we introduce the following reference equation:
\beg{align}\label{FrozeEQ}
\d\tilde{X}_t=b(t,\tilde{X}_t,\mu_t)\d t+\sigma(t,\mu_t)\d B_t^H,\ \ 0\leq t\leq T
\end{align}
with initial value $\tilde{X}_0=y\in\R^d$.
With the help of the Girsanov theorem for the fractional Brownian motion, we can obtain the following perturbation result, whose proof is postponed to the Appendix.
\beg{lem}\label{Le(Meas)}
Suppose that \textsc{\textbf{(H1)}} holds.
Then for any $\mu_\cdot\in C([0,T]; \sP_2(\R^d))$, there is a measurable map $\G_\mu: C([0,T]; \R^d)\rightarrow C([0,T]; \R^d)$ such that
\begin{align*}
\tilde{X}_\cdot=\mathcal{G}_\mu(B_\cdot^H).
\end{align*}
Moreover, for each $h\in\mathcal{A}_M$, define
\begin{align*}
\tilde{X}_\cdot^h:=\mathcal{G}_\mu\left(B_\cdot^H+(R_Hh)(\cdot)\right),
\end{align*}
then $\tilde{X}^h$ satisfies the following equation
\begin{align}\label{1-Le(Meas)}
\tilde{X}_t^h=&y+\int_0^tb(s, \tilde{X}_s^h, \mu_s)\d s+\int_0^t\si(s, \mu_s)\d(R_Hh)(s)\cr
&+\int_0^t\si(s, \mu_s)\d B_s^H,\ \ t\in[0,T],\ \ \P\textit{-}a.s..
\end{align}
\end{lem}
From the above lemma, it easily follows the following result.
\beg{lem}\label{Le(Meas2)}
Suppose that $y=x$ and $\mu_t=\sL_{X_t}, t\in[0,T]$ for equation \eqref{FrozeEQ} and \textsc{\textbf{(H1)}} holds.
Then the solution $X$ of equation \eqref{DDsde} satisfies $X_\cdot=\mathcal{G}_{\sL_{X}}(B_\cdot^H)$, where $\mathcal{G}_{\sL_{X}}$ is given in Lemma \ref{Le(Meas)} with $\mu=\sL_{X}$.
Moreover, for any $h\in\mathcal{A}_M$, let
\begin{align*}
X_\cdot^h=\mathcal{G}_{\sL_X}\left(B_\cdot^H+(R_Hh)(\cdot)\right),
\end{align*}
then $X^h$ satisfies the following equation
\begin{align*}
X_t^h=&x+\int_0^tb(s, X_s^h, \sL_{X_s})\d s+\int_0^t\sigma(s, \sL_{X_s})\d(R_Hh)(s)\cr
&+\int_0^t\sigma(s,\sL_{X_s})\d B_s^H,\ \ t\in[0,T],\ \ \P\textit{-}a.s..
\end{align*}
\end{lem}

In this article, our main objective is to study asymptotic behaviors for DDSDEs driven by fractional Brownian motions.
More precisely, we will investigate the large and moderate deviation principles and the central limit theorem of the solution to  the following DDSDE with small fractional noise: for any $\ep>0$,
\begin{align}\label{Ap-DDsde}
\d X_t^\epsilon =b(t,X_t^\ep,\mathscr{L}_{X_t^\ep})\d t+\ep^H\sigma(t,\mathscr{L}_{X_t^\ep})\d B_t^H,\ \ X_0^\ep =x.
\end{align}
To this end, we first present preliminary fact and result.
According to Lemma \ref{Le(Meas2)}, there exists a measurable map $\mathcal{G}^\epsilon:= \mathcal{G}_{\mathscr{L}_{X^\epsilon}}$ such that $X_\cdot^\epsilon=\mathcal{G}^\epsilon(\ep^H B_\cdot^H).$
Furthermore, for every $h^\ep\in\mathcal{A}_M$, let
\begin{align}\label{Per2-DDsde}
X_\cdot^{\ep, h^\ep}:=\mathcal{G}^\ep\left(\ep^H B_\cdot^H+(R_Hh^\ep)(\cdot)\right),
\end{align}
then $X^{\ep, h^\ep}$ satisfies the following equation
\begin{align}\label{Per-DDsde}
X_t^{\ep, h^\ep}=&x+\int_0^tb(s, X_s^{\ep, h^\ep}, \sL_{X_s^\ep})\d s+\int_0^t\sigma(s,\sL_{X_s^\ep})\d(R_Hh^\ep)(s)\cr
&+\ep^H\int_0^t\sigma(s, \sL_{X_s^\ep})\d B_s^H,\ \ t\in[0,T],\ \ \P\textit{-}a.s..
\end{align}
Besides, again due to \cite[Theorem 3.1]{FHSY}, we have the following result.

\beg{prp}\label{ODE}
Suppose that \textsc{\textbf{(H1)}} holds.
Then there exists a unique function $\{X_t^0\}_{t\in[0,T]}$ such that
\begin{itemize}
\item[(i)]$X^0\in C([0,T];\R ^d)$,
\item[(ii)]$X^0$ satisfies the following deterministic equation
\begin{align}\label{LimSDE}
X_t^0=x+\int_0^t b(s, X_s^0, \sL_{X_s^0})\d s, \ \ t\in[0,T].
\end{align}
\end{itemize}
\end{prp}
It is easy to see that $X^0$ is deterministic and $\sL_{X_s^0}=\de_{X_s^0}$.
In the sequel, we always use $X^0$ to denote the unique solution of equation \eqref{LimSDE}.

\subsection{Large deviation principle (LDP)}

In this part, we aim to establish the LDP for equation \eqref{Ap-DDsde} as $\ep\ra0$.
We first recall some definitions of the theory of LDP.
Let $\mathscr{E}$ be a Polish space with the Borel $\sigma$-field $\mathcal{B}(\mathscr{E})$.

\beg{defn}\label{De(rate)}
(Rate function)
A function $I$ is called a rate function, if for each constant $M<\infty$, the level set $\{x\in\mathscr{E}: I(x)\leq M\}$ is a compact subset of $\mathscr{E}$.
\end{defn}

\beg{defn}\label{De(ldp)}
(Large deviation principle)
Let $I$ be a rate function on $\mathscr{E}$.
Given a collection $\{\ell(\ep)\}_{\ep>0}$ of positive reals,
a family $\{\mathbb{X}^\ep\}_{\ep>0}$ of $\mathscr{E}$-valued random variables is said to be satisfied a LDP on $\mathscr{E}$ with speed $\ell(\ep)$ and rate function $I$
if the following two conditions hold:
\begin{itemize}
\item[(i)](Upper bound) For each closed subset $F\subset\mathscr{E}$,
\begin{align*}
\limsup_{\ep\to 0}\ell(\ep)\log{\mathbb{P}(\mathbb{X}^\ep\in F)} \leq -\inf_{x\in F}I(x).
\end{align*}

\item[(ii)](Lower bound) For each open subset $G\subset\mathscr{E}$,
\begin{align*}
\liminf_{\ep\to 0}\ell(\ep)\log{\mathbb{P}(\mathbb{X}^\ep\in G)} \geq -\inf_{x\in G}I(x).
\end{align*}
\end{itemize}
\end{defn}
By \cite[Theorems 1.5 and 1.8]{BD19} (see also \cite[Theorems 1.2.1 and 1.2.3]{DE11}),  the large deviation principle is equivalent to the following so-called Laplace principle.

\beg{defn}\label{De(lp)}(Laplace principle)
Let $I$ be a rate function on $\mathscr{E}$. Given a collection $\{\ell(\ep)\}_{\ep>0}$ of positive reals,
a family $\{\mathbb{X}^\ep\}_{\ep>0}$ of $\mathscr{E}$-valued random variables is said to be satisfied the Laplace principle upper bound (respectively, lower bound) on $\mathscr{E}$ with speed $\ell(\ep)$ and rate function $I$
if for all $\varrho\in C_b(\mathscr{E})$,
\begin{align}\label{De(lp)1}
\limsup_{\ep\to 0}-\ell(\ep)\log\E\left[\exp\left(-\ff {\varrho(\mathbb{X}^\ep)}{\ell(\ep)} \right)\right]\leq\inf_{x\in\mathscr{E}}\{\varrho(x)+I(x)\},
\end{align}
(respectively,
\begin{align}\label{De(lp)2}
\liminf_{\ep\to 0}-\ell(\ep)\log\E\left[\exp\left(-\ff {\varrho(\mathbb{X}^\ep)}{\ell(\ep)} \right)\right]\geq\inf_{x\in\mathscr{E}}\{\varrho(x)+I(x)\}).
\end{align}
The Laplace principle is said to be held for $\{\mathbb{X}^\ep\}$  with speed $\ell(\ep)$ and rate function $I$ if both the Laplace upper and lower bounds hold.
\end{defn}

For any $\ep>0$, let $\G^\ep: C([0,T]; \R^d)\ra\mathscr{E}$ be a measurable map (with a slight abuse of notation $\G^\ep$).
Next, we give the following sufficient condition for the Laplace principle (equivalently, the LDP) of $\mathbb{X}^\ep=\G^\ep(\ep^H B_\cdot^H)$ as $\ep\ra0$,
which is a fractional version of \cite[Theorem 4.2]{BDM11} or \cite[Theorem 4.3]{LSZZ22}.

\begin{enumerate}

\item[\textsc{\textbf{(A0)}}] There exists a measurable map  $\mathcal{G}^0: I_{0+}^{H+1/2}(L^2([0,T],\R^d))\ra \mathscr{E}$ such that the following two conditions hold.

\item[(i)] Let $\{h^\ep:\ep>0\}\subset\mathcal{A}_M$ for any $M\in(0,\infty)$. If $h^\ep$ converges to $h$ in distribution as $S_M$-valued random elements, then
\begin{align*}
\G^\ep(\ep^HB^H_\cdot+\ep^H/\ell^{\ff 1 2}(\ep)(R_Hh^\ep)(\cdot))\ra \G^0(R_Hh)
\end{align*}
in law as $\ep\ra0$, where $\{\ell^\ep\}_{\ep>0}$ are positive reals.

\item[(ii)]  For each $M\in(0,\infty)$, the set $\{\G^0(R_Hh):h\in S_M\}$ is a compact subset of $\mathscr{E}$.

\end{enumerate}

\beg{prp}\label{Suf1(LDP)}
If  $\mathbb{X}_\cdot^\ep=\G^\ep(\ep^H B_\cdot^H)$ and \textsc{\textbf{(A0)}} holds, then the family $\{\mathbb{X}^\ep:\ep>0\}$ satisfies the Laplace principle (hence the LDP) on $\mathscr{E}$ with speed $\ell(\ep)$
and the rate function $I$ given by
\begin{align}\label{1-Suf1(LDP)}
I(f)=\inf_{\{h\in\H:f=\G^0(R_Hh)\}}\left\{\ff 1 2 \|h\|_{\H}^2\right\},\ \ f\in \mathscr{E}.
\end{align}
Here we follow the convention that the infimum over an empty set is $+\infty$.
\end{prp}

For the sake of conciseness, we defer the proof to the Appendix.
Below is a convenient and sufficient condition for verifying \textsc{\textbf{(A0)}} in Proposition \ref{Suf1(LDP)}.
The proof is pretty similar to that of \cite[Theorem 3.2]{MSZ21} and we omit it here.

\begin{enumerate}

\item[\textsc{\textbf{(A1)}}] There exists a measurable map  $\G^0: I_{0+}^{H+1/2}(L^2([0,T],\R^d))\ra\mathscr{E}$ for which the following two conditions hold.

\item[(i)] Let $\{h^\ep:\ep>0\}\subset\mathcal{A}_M$ for any $M\in(0,\infty)$. For each $\delta>0$,
\begin{align*}
\lim_{\ep\ra 0}\P\left(d\left(\G^\ep(\ep^H B_\cdot^H+\ep^H/\ell^{\ff 1 2}(\ep)(R_Hh^\ep)(\cdot)),\G^0((R_Hh^\ep)(\cdot))\right)>\delta\right)=0,
\end{align*}
where $d(\cdot,\cdot)$ stands for the metric on $ \mathscr{E}, \{\ell^\ep\}_{\ep>0}$ are positive reals.

\item[(ii)]  Let $\{h^n:n\in\mathbb{N}\}\subset S_M$ for any $M\in(0,\infty)$. If $h^n$ converges to some element $h$ in $S_M$ as $n\ra\infty$, then $\G^0(R_Hh^n)$ converges to $\G^0(R_Hh)$ in $\mathscr{E}$.

\end{enumerate}

\beg{prp}\label{Suf2(LDP)}
If $\mathbb{X}_\cdot^\epsilon=\mathcal{G}^\ep(\ep^H B_\cdot^H)$ and \textsc{\textbf{(A1)}} holds, then the family $\{\mathbb{X}^\ep:\ep>0\}$ satisfies the Laplace principle (hence the LDP) on $\mathscr{E}$ with speed $\ell(\ep)$
and the rate function $I$ given by \eqref{1-Suf1(LDP)}.
\end{prp}

Before moving to the LDP for equation \eqref{Ap-DDsde}, we consider the following skeleton equation
\begin{align}\label{SkEq}
\Upsilon_t^h=x+\int_0^tb(s,\Upsilon_s^h,\sL_{X_s^0})\d s+\int_0^t\sigma(s,\sL_{X_s^0})\d(R_H h)(s),\ \ t\in[0,T],
\end{align}
where $h\in\H$ and $X^0$ is given in \eqref{LimSDE}.
Here and in the next part, we further assume that $\si$ is H\"{o}lder continuous of order belonging to $(1-H,1]$ with respect to the time variable.
Then,  $\int_0^t\sigma(s,\sL_{X_s^0})\d(R_H h)(s)$ in \eqref{SkEq} (or \eqref{Sk(mdp)} below) is well-defined as a Riemann-Stieltjes integral because of \eqref{HORH} and Proposition \ref{ODE}.
For any $t\in[0,T]$ and $x\in\R^d$, set
\begin{align*}
\bar{b}(t,x):=b(t,x,\sL_{X_t^0}),\ \ \ \ \ \bar{\si}(t):=\sigma(t,\sL_{X_t^0}).
\end{align*}
It is easy to see that $\bar{b}(t,x)$ and $\bar{\si}(t)$ independent of the measure satisfy \textsc{\textbf{(H1)}},
which ensure that equation \eqref{SkEq} has a unique solution due to \cite[Theorem 5.1]{NR02}.
As a consequence, we can define a map as follows
\begin{align}\label{RateF}
\G^0: I_{0+}^{H+1/2}(L^2([0,T],\R^d))\ni R_Hh\mapsto \Upsilon^h\in C([0,T]; \R^d),
\end{align}
where $h\in\H$ and $\Upsilon^h$ is the unique solution of equation \eqref{SkEq}.
Besides, one can see that as $\ep\ra0$, equation \eqref{Ap-DDsde} reduces to \eqref{LimSDE} and then $\sL_{X_\cdot^\ep}$ goes to $\sL_{X_\cdot^0}$.
So, in this sense $\sL_{X_\cdot^0}$, rather than $\sL_{\Upsilon_\cdot^h}$, appearing in equation \eqref{SkEq} to define the rate function of Theorem \ref{Th(ldp)} below, seems to reasonable.

Our main result in this part reads as follows.

\beg{thm}\label{Th(ldp)}
Assume that \textsc{\textbf{(H1)}} holds. For each $\ep>0$, let $X^\ep=\{X_t^\ep\}_{t\in[0,T]}$ be the solution to  equation \eqref{Ap-DDsde}.
Then the family $\{X^\ep:\ep>0\}$ satisfies a LDP on $C([0,T]; \R^d)$ with speed $\ep^{2H}$ and the  rate function $I$ given by \eqref{1-Suf1(LDP)},
where $\G^0$ is defined in \eqref{RateF}.
\end{thm}

\subsection{Moderate deviation principle (MDP)}

In this part, we shall investigate the MDP for the equation \eqref{Ap-DDsde} as $\ep\ra0$.
The moderate deviations problem for $\{X^\ep:\ep>0\}$ is to study the asymptotics of
\begin{align*}
\ff 1{\kappa^2(\ep)}\log\P(Y^\ep\in\cdot),
\end{align*}
where $\kappa(\ep)\ra\infty, \ep^H\kappa(\ep)\ra0$ as $\ep\ra0$ and
\begin{align}\label{Meq-0}
Y^\ep:=\ff {X^\ep-X^0}{\ep^H\kappa(\ep)}.
\end{align}
Recall the equations \eqref{Ap-DDsde} and \eqref{LimSDE}, from which we easily deduce that $Y^\ep$ satisfies
\begin{align}\label{Meq-1}
Y^\ep_t=&\ff 1{\ep^H\kappa(\ep)}\int_0^t(b(t,X^0_s+\ep^H\kappa(\ep)Y_s^\ep,\mathscr{L}_{X_s^\ep})-b(s, X_s^0, \mathscr{L}_{X_s^0}))\d s\cr
&+\ff 1 {\kappa(\ep)}\int_0^t\si(s,\mathscr{L}_{X_s^\ep})\d B_s^H,\ \ t\in[0,T].
\end{align}

Next, we put
\begin{align*}
\widetilde{\G}^\ep(\cdot):=\ff {\mathcal{G}^\ep(\cdot)-X^0}{\ep^H\kappa(\ep)},
\end{align*}
which is a map from $C([0,T];\R^d)$ to  $C([0,T];\R^d)$ such that $Y^\ep=\widetilde{\G}^\ep(\ep^H B_\cdot^H)$
due to the definition of $\G^\epsilon$ and the relation $X_\cdot^\ep=\G^\ep(\ep^H B_\cdot^H)$.
Moreover, for any $h^\ep\in\mathcal{A}_M$, let
\begin{align}\label{Add(Meq-2)}
Y_\cdot^{\ep, h^\ep}=\widetilde{\G}^\ep\left(\ep^H B_\cdot^H+\ep^H\kappa(\ep)(R_Hh^\ep)(\cdot)\right),
\end{align}
then $Y^{\ep, h^\ep}$ solves the following equation
\begin{align}\label{Meq-2}
Y_t^{\ep, h^\ep}=&\ff 1{\ep^H\kappa(\ep)}\int_0^t\left(b(s, X^0_s+\ep^H\kappa(\ep)Y_s^{\ep, h^\ep}, \sL_{X_s^\ep})- b(s, X_s^0, \mathscr{L}_{X_s^0})\right)\d s\cr
&+\int_0^t\si(s,\sL_{X_s^\ep})\d(R_Hh^\ep)(s)
+\ff 1{\kappa(\ep)}\int_0^t\si(s, \sL_{X_s^\ep})\d B_s^H,\ \ t\in[0,T],\ \ \mathbb{P}\textit{-}a.s..
\end{align}

To obtain the MDP, in additional to \textsc{\textbf{(H1)}}, we also need the following assumption.
\begin{enumerate}

\item[\textsc{\textbf{(H2)}}] The derivative $\nabla b(t,\cdot,\mu)(x)$ exists and there is a non-decreasing function $\widetilde{K}(t)$ such that for any $t\in[0,T], x,y\in\R^d, \mu\in\mathcal{P}_\th(\R^d)$,
\begin{align*}
\left\|\nabla b(t,\cdot, \mu)(x)-\nabla b(t,\cdot, \mu )(y) \right\|\leq \widetilde{K}(t)(|x-y|).
\end{align*}

\end{enumerate}

Now, for each $h\in\H$, we introduce the following equation
\begin{align}\label{Sk(mdp)}
\Xi_t^h=\int_0^t\nabla_{\Xi_s^h} b(s,\cdot,\sL_{X_s^0})(X_s^0)\d s+\int_0^t\sigma(s,\sL_{X_s^0})\d(R_H h)(s),\ \ t\in[0,T],
\end{align}
which is used to give the rate function of Theorem \ref{Th(mdp)} below.
Similar to the above part, under the time H\"{o}lder continuity of $\si$ with order belonging to $(1-H,1]$,
the equation \eqref{Sk(mdp)} admits a unique solution, in which the last integral is also regarded as a Riemann-Stieltjes integral.
Therefore, this allows us to define a map as follows
\begin{align}\label{RateF(mdp)}
\widetilde{\G}^0: I_{0+}^{H+1/2}(L^2([0,T],\R^d))\ni R_Hh\mapsto \Xi^h\in C([0,T]; \R^d),
\end{align}
where $h\in\H$ and $\Xi^h$ is the unique solution of the equation \eqref{Sk(mdp)}.

We state our main result of this part as follows.

\beg{thm}\label{Th(mdp)}
Assume that \textsc{\textbf{(H1)}} and \textsc{\textbf{(H2)}} hold. For each $\ep>0$, let $Y^\ep=\{Y_t^\ep\}_{t\in[0,T]}$ be defined in \eqref{Meq-0}.
Then the family $\{Y^\ep:\ep>0\}$ satisfies a LDP on $C([0,T]; \R^d)$ with speed $\kappa^{-2}(\ep)$ and the rate function $I$ given by \begin{align}\label{1-Th(mdp)}
I(f)=\inf_{\{h\in\H:f=\widetilde{\G}^0(R_Hh)\}}\left\{\ff 1 2 \|h\|_{\H}^2\right\},\ \ f\in C([0,T]; \R^d),
\end{align}
where $\widetilde{\G}^0$ is defined in \eqref{RateF(mdp)}.
Here we use the convention that the infimum over an empty set is $+\infty$.
\end{thm}

\subsection{Central limit theorem (CLT)}

This part is devoted to studying the CLT for equation \eqref{Ap-DDsde}.
More precisely, we shall show that $\ff {X^\ep-X^0}{\ep^H}$ converges to a stochastic process in the $p$-moment sense as $\ep\ra0$.
We would like to mention that the limit process is a solution to some linear equation which involves the Lions derivative of the coefficient $b$.
To this end, we will impose the following conditions on $b$ and $\si$.

\begin{enumerate}

\item[\textsc{\textbf{(H3)}}] For every $t\in[0,T]$, $b(t,\cdot,\cdot)\in C^{1,(1,0)}(\R^d\times\sP_2(\R^d))$, and there exists a non-decreasing function $\bar{K}(t)$ such that

\item[(i)] for any $t\in[0,T],\ x,y\in\R^d,\ \mu,\nu\in\sP_2(\R^d)$,
\begin{align*}
\|\nabla b(t,\cdot,\mu)(x)\|+|D^Lb(t,x,\cdot)(\mu)(y)|\le \bar{K}(t),\ \ \|\si(t,\mu)-\si(t,\nu)\|\leq\bar{K}(t)\W_\th(\mu,\nu),
\end{align*}
and $|b(t,0,\de_0)|+\|\si(t,\de_0)\|\leq \bar{K}(t)$.

\item[(ii)]  for any $t\in[0,T],\ x,y,z_1,z_2\in\R^d,\ \mu,\nu\in\sP_2(\R^d)$,
\begin{align*}
&\|\nabla b(t,\cdot,\mu)(x)-\nabla b(t,\cdot,\nu)(y)\|+|D^Lb(t,x,\cdot)(\mu)(z_1)-D^Lb(t,y,\cdot)(\nu)(z_2)|\cr
&\le \bar{K}(t)(|x-y|+|z_1-z_2|+\W_\th(\mu,\nu)).
\end{align*}

\end{enumerate}
Observe that according to the fundamental theorem for Bochner integral (see, for instance, \cite[Proposition A.2.3]{LWbook}) and the definitions of $L$-derivative and the Wasserstein distance,
\textsc{\textbf{(H3)}}(i) implies
\begin{align*}
|b(t,x,\mu)-b(t,y,\nu)|\le \bar{K}(t)(|x-y|+\W_\th(\mu,\nu)),\ \  t\in[0,T],\ x,y\in\R^d,\ \mu,\nu\in\sP_2(\R^d).
\end{align*}
Then, it follows from \cite[Theorem 3.1]{FHSY} that \eqref{Ap-DDsde} has a unique solution.

Next, we give an example of the function $b$ such that \textsc{\textbf{(H3)}} is satisfied.

\beg{exa}\label{Ex-clt}
Assume  $f:[0,T]\times\R^d\ra\R^d$ and $\varphi:\R^d\ra\R^d$ are all twice continuously differentiable mappings with bounded derivatives.
Suppose that
\beg{align*}
b(t,x,\mu)=f\left(t,x+\int_{\R^d}\varphi(u)\mu(\d u)\right).
\end{align*}
Then, it is easy to see that
\beg{align*}
&\na b(t,\cdot,\mu)(x)=\na f\left(t,x+\int_{\R^d}\varphi(u)\mu(\d u)\right),\cr
&D^Lb(t,x,\cdot)(\mu)(y)=\na f\left(t,x+\int_{\R^d}\varphi(u)\mu(\d u)\right)\na\varphi(y).
\end{align*}
In particular, when $\varphi=0, b$ reduces to a function with no dependence on the measure.
By a direct calculation, it is readily checked that the function $b$ above satisfies  \textsc{\textbf{(H3)}}.
More  examples can be found in our recent work \cite[Example 4.6]{FHSY}.
\end{exa}

Our main result in this part is stated in the following theorem.

\beg{thm}\label{Th(clt)}
Assume that \textsc{\textbf{(H3)}} holds, then for any $p\geq\th$ and $p>1/H$,
\begin{align*}
\E\left (\sup_{0\leq t\leq T}\left|\frac{X_t^\epsilon-X_t^0}{\epsilon^H}-Z_t \right|^p \right) \leq
C_{T,p,H}\ep^{pH}\bigg(1+\sup_{t\in[0,T]}|X_t^0|^{2p}\bigg),\ \  \ep\in(0,\ep_0],
\end{align*}
where $Z_t$ satisfies
\begin{align}\label{1-Th(clt)}
Z_t=&\int_0^t\nabla_{Z_s}b(s,\cdot,\sL_{X_s^0})(X_s^0)\d s+\int_0^t\left(\E\langle D^Lb(s,u,\cdot)(\sL_{X_s^0})(X_s^0), Z_s \rangle\right)|_{u=X_s^0}\d s\cr
&+\int_0^t\sigma(s,\sL_{X_s^0})\d B_s^H,\ \ t\in[0,T],
\end{align}
and $\ep_0>0$ is a constant appeared in Lemma \ref{Dif(TE)} below.
\end{thm}


\section{Proofs of the main results}

To prove the main results, we need the following useful lemma, which presents a maximal inequality for $\int_0^t\sigma (s,\mu_s)\d B_s^H$.
Though its proof is identical to step 1 in the proof of \cite[Theorem 3.1]{FHSY}, yet we give its proof in the Appendix for the convenience of the reader.

\beg{lem}\label{MomEs}
Suppose that $\si$ satisfies \textsc{\textbf{(H1)}} and $\mu\in C([0,T]; \sP_p(\R^d))$ with $p\geq\th$ and $p>1/H$.
Then there is a constant $C_{T,p,H}>0$ such that
\begin{align}\label{MomEs-1}
\E\left(\sup_{t\in[0,T]}\left|\int_0^t\sigma (s,\mu_s)\d B_s^H \right|^p\right)\leq C_{T,p,H}\int _0^T\|\sigma (s,\mu_s)\|^p\d s.
\end{align}
\end{lem}

\subsection{Proof of Theorem \ref{Th(ldp)}}

According to Proposition \ref{Suf2(LDP)}, to prove Theorem \ref{Th(ldp)}, it is enough to check that \textsc{\textbf{(A1)}} holds
with $\G^\ep, \G^0$ and $\ell(\ep)$ given by \eqref{Per2-DDsde}, \eqref{RateF} and $\ep^{2H}$, respectively.
The verification of \textsc{\textbf{(A1)}}(i) and \textsc{\textbf{(A1)}}(ii) will be shown respectively in Proposition \ref{PrA12} and Proposition \ref{PrA11} below.

Before we give Proposition \ref{PrA12}, we first present the following lemma which characterizes the difference between $X^\epsilon$ and $X^0$.

\beg{lem}\label{Dif(TE)}
Suppose that \textsc{\textbf{(H1)}} holds.
Then for any $p\geq\th$ and $p>1/H$, there exists a constant $\ep_0>0$ such that for every $\ep\in(0,\ep_0]$,
\begin{align*}
\E\bigg(\sup_{t\in[0,T]}\left|X_t^\epsilon-X_t^0\right|^p\bigg)\leq C_{T,p,H}\ep^{pH}\bigg(1+\sup_{t\in[0,T]}|X_t^0|^p\bigg),
\end{align*}
\end{lem}
\beg{proof}
By \eqref{Ap-DDsde}, \eqref{LimSDE} and  \textsc{\textbf{(H1)}}, we have for any $t\in[0,T]$,
\begin{align}\label{1-PfDif(TE)}
&\sup_{s\in[0,t]}\left|X_s^\epsilon-X_s^0\right|^p\cr
\leq& 2^{p-1}\sup_{s\in[0,t]}\left|\int_0^s\left(b(r, X_r^\epsilon, \sL_{X_r^\epsilon})-b(r, X_r^0, \sL_{ X_r^0})\right)\d r\right|^p\cr
&+2^{p-1}\ep^{pH}\sup_{s\in[0,t]}\left|\int_0^s\sigma(r,\sL_{X_r^\epsilon})\d B_r^H\right|^p\cr
\leq&(4T)^{p-1}K^p(T)\int_0^t\left(\left|X_r^\epsilon-X_r^0\right|^p+\W_\th^p(\sL_{X_r^\epsilon},\sL_{X_r^0})\right)\d r\cr
&+2^{p-1}\ep^{pH}\sup_{s\in[0,t]}\left|\int_0^s\sigma(r,\sL_{X_r^\epsilon})\d B_r^H\right|^p.
\end{align}
Note that by Lemma \ref{MomEs} and  \textsc{\textbf{(H1)}}, we get for any $p>1/H$,
\begin{align*}
&\E\left(\sup_{s\in[0,t]}\left|\int_0^s\sigma(r,\sL_{X_r^\epsilon})\d B_r^H\right|^p\right)\cr
\leq& C_{T,p,H}\int_0^t\|\si(r,\sL_{X_r^\ep})\|^p\d r\cr
\leq& 3^{p-1}C_{T,p,H}\int_0^t\|\si(r,\sL_{X_r^\ep})-\si(r,\sL_{X_r^0})\|^p\d r
+3^{p-1}C_{T,p,H}\int_0^t\|\si(r,\sL_{X_r^0})-\si(r,\delta_0)\|^p\d r\cr
&+3^{p-1}C_{T,p,H}\int_0^t\|\si(r,\delta_0)\|^p\d r\cr
\leq&3^{p-1}C_{T,p,H}K^p(T)\int_0^t\E\left(\sup_{u\in[0,r]}|X_u^\ep-X_u^0|^p\right)\d r+3^{p-1}C_{T,p,H}TK^p(T)\bigg(1+\sup_{t\in[0,T]}|X_t^0|^p\bigg).
\end{align*}
Then, taking the expectation on both sides of \eqref{1-PfDif(TE)} and using the Gronwall inequality, we derive that for any $p\geq\th$ and $p>1/H$,
\begin{align*}
\E\bigg(\sup_{t\in[0,T]}\left|X_t^\epsilon-X_t^0\right|^p\bigg)\leq C_{T,p,H}\e^{C_{T,p,H}(1+\ep^{pH})}\ep^{pH}\bigg(1+\sup_{t\in[0,T]}|X_t^0|^p\bigg).
\end{align*}
Hence, there exists a constant $\ep_0>0$ such that for every $\ep\in(0,\ep_0]$,
\begin{align*}
\E\bigg(\sup_{t\in[0,T]}\left|X_t^\epsilon-X_t^0\right|^p\bigg)\leq C_{T,p,H}\ep^{pH}\bigg(1+\sup_{t\in[0,T]}|X_t^0|^p\bigg),
\end{align*}
which completes the proof.
\end{proof}

\beg{prp}\label{PrA12}
Suppose that  \textsc{\textbf{(H1)}} holds and let $\{h^\ep:\ep>0\}\subset\mathcal{A}_M$ for any $M\in(0,\infty)$.
Then, for any $\delta>0$,
\begin{align*}
\lim_{\ep\ra 0}\P\left(\|X_\cdot^{\ep, h^\ep}-\G^0((R_Hh^\ep)(\cdot))\|_\infty>\delta\right)=0,
\end{align*}
where $\|\cdot\|_\infty$ is the uniform norm on $C([0,T];\R^d)$.
\end{prp}
\beg{proof}
For each fixed $\ep>0$, by \eqref{SkEq}-\eqref{RateF} with $h$ replaced by $h^\ep$ and \eqref{Per-DDsde} we get
\begin{align*}
&X_t^{\ep, h^\ep}-\G^0((R_Hh^\ep)(\cdot))(t)=X_t^{\ep, h^\ep}-\Upsilon_t^{h^\ep}\cr
=&\int_0^t\left(b(s, X_s^{\ep, h^\ep}, \sL_{X_s^\ep})-b(s,\Upsilon_s^{h^\ep},\sL_{X_s^0})\right)\d s\cr
&\ \ +\int_0^t\left(\sigma(s,\sL_{X_s^\ep})-\sigma(s,\sL_{X_s^0})\right)\d(R_Hh^\ep)(s)
 +\ep^H\int_0^t\sigma(s, \sL_{X_s^\ep})\d B_s^H,\ \ t\in[0,T].
\end{align*}
Then, it follows that
\begin{align}\label{1-PfPrA12}
|X_t^{\ep, h^\ep}-\Upsilon_t^{h^\ep}|^2\leq&
3\left|\int_0^t\left(b(s, X_s^{\ep, h^\ep}, \sL_{X_s^\ep})-b(s,\Upsilon_s^{h^\ep},\sL_{X_s^0})\right)\d s\right|^2\cr
&+3\left|\int_0^t\left(\sigma(s,\sL_{X_s^\ep})-\sigma(s,\sL_{X_s^0})\right)\d(R_Hh^\ep)(s)\right|^2\cr
&+3\ep^{2H}\left|\int_0^t\sigma(s, \sL_{X_s^\ep})\d B_s^H\right|^2\cr
=:&J_1(t)+J_2(t)+J_3(t).
\end{align}
Using the H\"{o}lder inequality and  \textsc{\textbf{(H1)}}, we arrive at
\begin{align}\label{2-PfPrA12}
J_1(t)\leq&3\left|\int_0^tK(s)\left(|X_s^{\ep, h^\ep}-\Upsilon_s^{h^\ep}|+\W_\th(\sL_{X_s^\ep},\sL_{X_s^0})\right)\d s\right|^2\cr
\leq&6TK^2(T)\left(\int_0^t|X_s^{\ep, h^\ep}-\Upsilon_s^{h^\ep}|^2\d s+T\E\bigg({\sup_{t\in[0,T]}}|X_t^\ep-X_t^0|^2\bigg)\right).
\end{align}
With the help of \eqref{RFRH},  \textsc{\textbf{(H1)}} and the Fubini theorem, we deduce
\begin{align}\label{3-PfPrA12}
J_2(t)=&3\left|\int_0^t\left(\sigma(s,\sL_{X_s^\ep})-\sigma(s,\sL_{X_s^0})\right)\int_0^s\ff {\partial K_H}{\partial s}(s,r)(K_H^*h^\ep)(r)\d r\d s\right|^2\cr
\leq&3C_H^2\left(\int_0^tK(s)\W_\th(\sL_{X_s^\ep},\sL_{X_s^0})\int_0^s\left(\ff s r\right)^{H-\ff 1 2}(s-r)^{H-\ff 3 2}|(K_H^*h^\ep)(r)|\d r\d s\right)^2\cr
\leq&3C_H^2K^2(T)\E\bigg({\sup_{t\in[0,T]}}|X_t^\ep-X_t^0|^2\bigg)\cr
&\times\left(\int_0^tr^{\ff 1 2-H}|(K_H^*h^\ep)(r)|\left(\int_r^ts^{H-\ff 1 2}(s-r)^{H-\ff 3 2}\d s\right)\d r\right)^2\cr
\leq&C(T,H)\E\bigg({\sup_{t\in[0,T]}}|X_t^\ep-X_t^0|^2\bigg)\int_0^T|(K_H^*h^\ep)(r)|^2\d r\cr
=&C(T,H)\E\bigg({\sup_{t\in[0,T]}}|X_t^\ep-X_t^0|^2\bigg)\|h^\ep\|_{\H}^2,
\end{align}
where $C(T,H):=\ff{3C_H^2T^{2H}K^2(T)}{2(1-H)(H-\ff 1 2)^2}$ and the last equality is due to the fact that $K_H^*$ is an isometry between $\H$ and $L^2([0,T],\R^d)$.\\
For the term $J_3(t)$, from Lemma \ref{MomEs} and  \textsc{\textbf{(H1)}} we have
\begin{align}\label{4-PfPrA12}
&\E\bigg(\sup_{t\in[0,T]}J_3(t)\bigg)\cr
\leq&3C_{T,H}\ep^{2H}\int _0^T\|\sigma (s,\sL_{X_s^\ep})\|^2\d s\cr
\leq&9C_{T,H}\ep^{2H}\bigg(\int _0^T\|\sigma (s,\sL_{X_s^\ep})-\sigma (s,\sL_{X_s^0})\|^2\d s+\int _0^T\|\sigma (s,\sL_{X_s^0})-\sigma (s,\delta_0)\|^2
\d s\cr
&\qquad\qquad\quad+\int_0^T\|\sigma (s,\delta_0)\|^2\d s\bigg)\cr
\leq&9C_{T,H}TK^2(T)\ep^{2H}\left(\E\bigg({\sup_{t\in[0,T]}}|X_t^\ep-X_t^0|^2\bigg)+1+\sup_{t\in[0,T]}|X_t^0|^2\right).
\end{align}
Then, substituting  \eqref{2-PfPrA12}-\eqref{4-PfPrA12} into \eqref{1-PfPrA12} yields
\begin{align*}
&\E\bigg(\sup_{t\in[0,T]}|X_t^{\ep, h^\ep}-\Upsilon_t^{h^\ep}|^2\bigg)\cr
\leq&
6TK^2(T)\int_0^T \sup_{r\in[0,s]}|X_r^{\ep, h^\ep}-\Upsilon_r^{h^\ep}|^2\d s\cr
&+\left[6(TK(T))^2+C(T,H)\|h^\ep\|_{\H}^2+9C_{T,H}TK^2(T)\ep^{2H}\right]\E\bigg({\sup_{t\in[0,T]}}|X_t^\ep-X_t^0|^2\bigg)\cr
&+9C_{T,H}TK^2(T)\ep^{2H}\bigg(1+\sup_{t\in[0,T]}|X_t^0|^2\bigg).
\end{align*}
Applying the Gronwall inequality and Lemma \ref{Dif(TE)} and using $h^\ep\in\mathcal{A}_M$,
we conclude that there is a constant $C_{T,H,M}$ such that
\begin{align*}
\E\bigg(\sup_{t\in[0,T]}|X_t^{\ep, h^\ep}-\Upsilon_t^{h^\ep}|^2\bigg)\leq C_{T,H,M}\ep^{2H}(1+\ep^{2H})\bigg(1+\sup_{t\in[0,T]}|X_t^0|^2\bigg).
\end{align*}
Then we get
\begin{align*}
\lim_{\ep\ra0}\E\bigg(\sup_{t\in[0,T]}|X_t^{\ep, h^\ep}-\Upsilon_t^{h^\ep}|^2\bigg)=0,
\end{align*}
which implies the desired assertion.
\end{proof}

To verify \textsc{\textbf{(A1)}}(ii), we need the following priori estimate for the solution $\Upsilon^h$ to the skeleton equation \eqref{SkEq}.

\beg{lem}\label{Est(SkE)}
Suppose that  \textsc{\textbf{(H1)}} holds. Then for any $M>0$,
\begin{align*}
\sup_{h\in S_M}\sup_{t\in[0,T]}|\Upsilon_t^h|^2\leq C_{T,H,M},
\end{align*}
where $C_{T,H,M}$ is a positive constant only depending on $T, H, M$.
\end{lem}
\beg{proof}
 According to the change-of-variables formula \cite[Theorem 4.3.1]{Zahle} and  \textsc{\textbf{(H1)}}, we have
\begin{align}\label{1-PfEst(SkE)}
|\Upsilon_t^h|^2&=|x|^2+2\int_0^t\langle \Upsilon_s^h, b(s,\Upsilon_s^h,\sL_{X_s^0})\rangle\d s+2\int_0^t\langle \Upsilon_s^h, \sigma(s,\sL_{X_s^0})\d(R_H h)(s)\rangle\cr
&=:|x|^2+I_1(t)+I_2(t).
\end{align}
From  \textsc{\textbf{(H1)}}, it follows that
\begin{align}\label{2-PfEst(SkE)}
I_1(t)\leq&2\int_0^t|\Upsilon_s^h|\cdot|b(s,\Upsilon_s^h,\sL_{X_s^0})-b(s,0,\delta_0)|\d s+2\int_0^t|\Upsilon_s^h|\cdot|b(s,0,\delta_0)|\d s\cr
\leq&2\int_0^tK(s)\left(|\Upsilon_s^h|^2+|\Upsilon_s^h|\cdot|X_s^0|\right)\d s+2\int_0^tK(s)|\Upsilon_s^h|\d s\cr
\leq&4K(t)\int_0^t|\Upsilon_s^h|^2\d s+K(t)\left(t+\int_0^t|X_s^0|^2\d s\right).
\end{align}
For the term $I_2(t)$, using \eqref{RFRH} and the fact that $K_H^*$ is an isometry between $\H$ and $L^2([0,T],\R^d)$, we have
\begin{align}\label{0-PfEst(SkE)}
I_2(t)
=&2\sum_{i,j=1}^d\int_0^T\Upsilon_{i,s}^h\sigma_{ij}(s,\sL_{X_s^0})\mathrm{I}_{[0,t]}(s)\int_0^s\ff {\partial K_H}{\partial s}(s,r)(K_H^*h)_j(r)\d r\d s\cr
=&2\sum_{i,j=1}^d\int_0^T\int_r^T\Upsilon_{i,s}^h\sigma_{ij}(s,\sL_{X_s^0})\mathrm{I}_{[0,t]}(s) \ff {\partial K_H}{\partial s}(s,r)\d s(K_H^*h)_j(r)\d r\cr
=&2\sum_{i,j=1}^d\int_0^T(K_H^*(\Upsilon_{i,\cdot}^h\sigma_{ij}(\cdot,\sL_{X_\cdot^0}))\mathrm{I}_{[0,t]})(r)(K_H^*h)_j(r)\d r\cr
=&2\int_0^T\langle K_H^*(\si^T(\cdot,\sL_{X_\cdot^0})\Upsilon^h\mathrm{I}_{[0,t]})(r),(K_H^*h)(r)\rangle\d r\cr
=&2\langle\si^T(\cdot,\sL_{X_\cdot^0})\Upsilon^h\mathrm{I}_{[0,t]},h\rangle_\H.
\end{align}
Here and in the sequel, $\si^T$ denotes the transpose matrix of $\si$.
Then, by \eqref{EsH} and  \textsc{\textbf{(H1)}} we get
\begin{align}\label{3-PfEst(SkE)}
I_2(t)
\leq&2\|\si^T(\cdot,\sL_{X_\cdot^0})\Upsilon^h\mathrm{I}_{[0,t]}\|_\H\cdot\|h\|_\H\cr
\leq&2(2H)^{1/2}T^{H-1/2}\|\si^T(\cdot,\sL_{X_\cdot^0})\Upsilon^h\mathrm{I}_{[0,t]}\|_{L^2}\cdot\|h\|_\H\cr
\leq&2(2H)^{1/2}T^{H-1/2}K(t)\left(1+\sup_{t\in[0,T]}|X_t^0|\right)\left(\int_0^t|\Upsilon_s^h|^2\d s\right)^{\ff 1 2}\|h\|_\H\cr
\leq&2HT^{2H-1}K(t)\left(1+\sup_{t\in[0,T]}|X_t^0|\right)\|h\|^2_\H
+K(t)\left(1+\sup_{t\in[0,T]}|X_t^0|\right)\int_0^t|\Upsilon_s^h|^2\d s.
\end{align}
Substituting \eqref{2-PfEst(SkE)} and \eqref{3-PfEst(SkE)} into \eqref{1-PfEst(SkE)} yields
\begin{align*}
\sup_{t\in[0,T]}|\Upsilon_t^h|^2
&\leq |x|^2+K(T)\left(T+\int_0^T|X_s^0|^2\d s\right)+2HT^{2H-1}K(T)\left(1+\sup_{t\in[0,T]}|X_t^0|\right)\|h\|^2_\H\cr
&\ \ +K(T)\left(5+\sup_{t\in[0,T]}|X_t^0|\right)\int_0^t|\Upsilon_s^h|^2\d s.
\end{align*}
Consequently, for any $h\in S_M$, the Gronwall inequality implies
\begin{align*}
\sup_{t\in[0,T]}|\Upsilon_t^h|^2&\leq\left[|x|^2+K(T)\left(T+\int_0^T|X_s^0|^2\d s\right)+2HT^{2H-1}K(T)\left(1+\sup_{t\in[0,T]}|X_t^0|\right)\|h\|^2_\H\right]\cr
&\ \ \times\exp\left\{TK(T)\left(5+\sup_{t\in[0,T]}|X_t^0|\right)\right\}\cr
&\leq C_{T,H,M},
\end{align*}
where the constant $C_{T,H,M}$ only depends on $T, H, M$.
The proof is now complete.
\end{proof}

\beg{prp}\label{PrA11}
Suppose that  \textsc{\textbf{(H1)}}  holds and let $\{h^n:n\in\mathbb{N}\}\subset\mathcal{S}_M$ for any $M\in(0,\infty)$
such that $h^n$ converges to element $h$ in $S_M$ as $n\ra\infty$.
Then
\begin{align*}
\lim_{n\ra\infty}\sup_{t\in[0,T]}|\G^0(R_Hh^n)(t)-\G^0(R_Hh)(t)|=0.
\end{align*}
\end{prp}
\beg{proof}
For each $n\geq1$, let $\Upsilon^{h^n}$ be the solution of equation \eqref{SkEq} with $h$ replaced by $h^n$.
By \eqref{RateF}, there hold $\G^0(R_Hh^n)=\Upsilon^{h^n}$ and $\G^0(R_Hh)=\Upsilon^h$.

We first prove that $\{\Upsilon^{h^n}\}_{n\geq1}$ is relatively compact in $C([0,T]; \R^d)$.
With the help of the Arzel\`{a}-Ascoli theorem, it is enough to show that $\{\Upsilon^{h^n}\}_{n\geq1}$ is uniformly bounded and equi-continuous in $C([0,T]; \R^d)$.
By Lemma \ref{Est(SkE)}, there exists a constant $C_{T,H,M}>0$ such that
\begin{align}\label{0-PfPrA11}
\sup_{n\geq1}\sup_{t\in[0,T]}|\Upsilon_t^{h^n}|\leq C_{T,H,M},
\end{align}
which means that $\{\Upsilon^{h^n}\}_{n\geq1}$ is uniformly bounded in $C([0,T]; \R^d)$.

Now, we focus on dealing with equi-continuous of $\{\Upsilon^{h^n}\}_{n\geq1}$ in $C([0,T]; \R^d)$.
By \eqref{SkEq}, we deduce that for $0\leq s<t\leq T$,
\begin{align}\label{1-PfPrA11}
\Upsilon_t^{h^n}-\Upsilon_s^{h^n}=\int_s^tb(r,\Upsilon_r^{h^n},\sL_{X_r^0})\d r+\int_s^t\sigma(r,\sL_{X_r^0})\d(R_H{h^n})(r).
\end{align}
In view of  \textsc{\textbf{(H1)}} and \eqref{0-PfPrA11}, we get
\begin{align}\label{2-PfPrA11}
\left|\int_s^tb(r,\Upsilon_r^{h^n},\sL_{X_r^0})\d r\right|
\leq&\int_s^t\left|b(r,\Upsilon_r^{h^n},\sL_{X_r^0})-b(r,0,\delta_0)\right|\d r+\int_s^t|b(r,0,\delta_0)|\d r\cr
\leq&\int_s^tK(r)\left(|\Upsilon_r^{h^n}|+|X_r^0|\right)\d r+\int_s^tK(r)\d r\cr
\leq& K(T)\bigg(1+C_{T,H,M}+\sup_{r\in[0,T]}|X_r^0|\bigg)(t-s).
\end{align}
By \eqref{RFRH},  \textsc{\textbf{(H1)}} and the Fubini theorem, we have
\begin{align}\label{3-PfPrA11}
&\left|\int_s^t\sigma(r,\sL_{X_r^0})\d(R_H{h^n})(r)\right|\cr
=&\left|\int_s^t\sigma(r,\sL_{X_r^0})\int_0^r\ff {\partial K_H}{\partial r}(r,u)(K_H^*h^n)(u)\d u\d r\right|\cr
\leq&C_HK(T)\bigg(1+\sup_{r\in[0,T]}|X_r^0|\bigg)\int_s^t\int_0^r\left(\ff r u\right)^{H-\ff 1 2}(r-u)^{H-\ff 3 2}|(K_H^*h^n)(u)|\d u\d r\cr
=&C_HK(T)\bigg(1+\sup_{r\in[0,T]}|X_r^0|\bigg)
\bigg[\int_0^su^{\ff 1 2-H}|(K_H^*h^n)(u)|\left(\int_s^tr^{H-\ff 1 2}(r-u)^{H-\ff 3 2}\d r\right)\d u\cr
&\qquad\qquad\qquad\qquad\qquad\quad \
+\int_s^tu^{\ff 1 2-H}|(K_H^*h^n)(u)|\left(\int_u^tr^{H-\ff 1 2}(r-u)^{H-\ff 3 2}\d r\right)\d u\bigg].
\end{align}
Using the H\"{o}lder inequality and the relation $\|K_H^*h^n\|_{L^2}=\|h^n\|_\H$, we arrive at
\begin{align}\label{4-PfPrA11}
&\int_0^su^{\ff 1 2-H}|(K_H^*h^n)(u)|\left(\int_s^tr^{H-\ff 1 2}(r-u)^{H-\ff 3 2}\d r\right)\d u\cr
\leq&T^{H-\ff 1 2}\int_s^t(r-s)^{H-\ff 3 2}\d r\int_0^su^{\ff 1 2-H}|(K_H^*h^n)(u)|\d u\cr
\leq&\ff{T^{\ff 1 2}}{(H-\ff 1 2)\sqrt{2(1-H)}}\left(\int_0^T|(K_H^*h^n)(u)|^2\d u\right)^{\ff 1 2}(t-s)^{H-\ff 1 2}\cr
=&\ff{T^{\ff 1 2}}{(H-\ff 1 2)\sqrt{2(1-H)}}\|h^n\|_\H(t-s)^{H-\ff 1 2}
\end{align}
and
\begin{align}\label{5-PfPrA11}
&\int_s^tu^{\ff 1 2-H}|(K_H^*h^n)(u)|\left(\int_u^tr^{H-\ff 1 2}(r-u)^{H-\ff 3 2}\d r\right)\d u\cr
\leq&\ff{T^{H-\ff 1 2}}{H-\ff 1 2}\int_s^tu^{\ff 1 2-H}(t-u)^{H-\ff 1 2}|(K_H^*h^n)(u)|\d u\cr
\leq&\ff{\sqrt{\mathcal{B}(2-2H,2H)}T^{H-\ff 1 2}}{H-\ff 1 2}\|h^n\|_\H\sqrt{t-s}.
\end{align}
Then, combining \eqref{1-PfPrA11}-\eqref{5-PfPrA11} together and using the fact that
$\|h^n\|_\H\leq\sqrt{2M}$ due to $h^n\in S_M$,
one can see that $\{\Upsilon^{h^n}\}_{n\geq1}$ is equi-continuous in $C([0,T]; \R^d)$.
Hence, we have shown that $\{\Upsilon^{h^n}\}_{n\geq1}$ is relatively compact in $C([0,T]; \R^d)$,
which implies that for any subsequence of $\{\Upsilon^{h^n}\}_{n\geq1}$, we can extract a further subsequence (not relabelled)
such that  $\Upsilon^{h^n}$ converges to some $\bar{\Upsilon}$ in $C([0,T]; \R^d)$.

We claim that $\bar{\Upsilon}=\Upsilon^h$.
Once this is shown, by a standard subsequential argument we can conclude that the full sequence  $\Upsilon^{h^n}$ converges to $\Upsilon^h$ in $C([0,T];\R^d)$, which is the desired assertion.

It remains to show the claim. By  \textsc{\textbf{(H1)}}, we first have for every $t\in[0,T]$,
\begin{align*}
&\left|\int_0^tb(s,\Upsilon_s^{h^n},\sL_{X_s^0})\d s-\int_0^tb(s,\bar{\Upsilon}_s,\sL_{X_s^0})\d s\right|\cr
\leq & \int_0^tK(s)\left|\Upsilon_s^{h^n}-\bar{\Upsilon}_s\right|\d s\cr
\leq & T K(T)\sup_{t\in[0,T]}\left|\Upsilon_t^{h^n}-\bar{\Upsilon}_t\right|\ra0, \ \ n\ra\infty.
\end{align*}
Then, for each $t\in[0,T]$, we get
\begin{align}\label{6-PfPrA11}
\lim_{n\ra\infty}\int_0^tb(s,\Upsilon_s^{h^n},\sL_{X_s^0})\d s=\int_0^tb(s,\bar{\Upsilon}_s,\sL_{X_s^0})\d s.
\end{align}
On the other hand, by \eqref{RFRH} and the Fubini theorem we derive
\begin{align}\label{7-PfPrA11}
&\int_0^t\sigma(s,\sL_{X_s^0})\d(R_H{h^n})(s)-\int_0^t\sigma(s,\sL_{X_s^0})\d(R_H{h})(s)\cr
=&\int_0^t\sigma(s,\sL_{X_s^0})\left(\int_0^s\ff {\partial K_H}{\partial s}(s,r)((K_H^*{h^n})(r)-(K_H^*h)(r))\d r\right)\d s\cr
=&C_H\int_0^T\left[\mathrm{I}_{[0,t]}(r)r^{\ff 1 2-H}\left(\int_r^t\sigma(s,\sL_{X_s^0})s^{H-\ff 1 2}(s-r)^{H-\ff 3 2}\d s\right)\right]\cr
&\qquad\quad\times[(K_H^*{h^n})(r)-(K_H^*h)(r)]\d r.
\end{align}
We set, for any unit vector $e\in\R^d$ and $t\in[0,T]$,
\begin{align*}
g_t(r):=\mathrm{I}_{[0,t]}(r)r^{\ff 1 2-H}\left(\int_r^t\sigma^T(s,\sL_{X_s^0})s^{H-\ff 1 2}(s-r)^{H-\ff 3 2}\d s\right)e,\ \ r\in[0,T].
\end{align*}
By  \textsc{\textbf{(H1)}} we have
\begin{align*}
|g_t(r)|\leq\ff{T^{2H}K(T)}{H+\ff 1 2}\bigg(1+\sup_{s\in[0,T]}\left|X_s^0\right|\bigg)r^{\ff 1 2-H},
\end{align*}
which implies that $g_t(\cdot)\in L^2([0,T],\R^d)$.
Then, taking into account of \eqref{7-PfPrA11} and the condition that $h^n$ converges to element $h$ in $S_M$ as $n\ra\infty$,
we deduce
\begin{align}\label{8-PfPrA11}
\lim_{n\ra\infty}\int_0^t\sigma(s,\sL_{X_s^0})\d(R_H{h^n})(s)=\int_0^t\sigma(s,\sL_{X_s^0})\d(R_H{h})(s).
\end{align}
Observe that $\Upsilon^{h^n}$ satisfies the following equation:
\begin{align*}
\Upsilon_t^{h^n}=x+\int_0^tb(s,\Upsilon_s^{h^n},\sL_{X_s^0})\d s+\int_0^t\sigma(s,\sL_{X_s^0})\d(R_H {h^n})(s),\ t\in[0,T],
\end{align*}
Letting $n$ goes to infinity and using \eqref{6-PfPrA11} and \eqref{8-PfPrA11}, one can see that $\bar{\Upsilon}$ solves \eqref{SkEq},
which yields $\bar{\Upsilon}=\Upsilon^h$ due to the uniqueness of solutions to \eqref{SkEq}.
This completes our proof.
\end{proof}

\subsection{Proof of Theorem \ref{Th(mdp)}}

In view of Proposition \ref{Suf2(LDP)}, to complete the proof of Theorem \ref{Th(mdp)}, it is sufficient to verify that  \textsc{\textbf{(A1)}} holds with $\G^0, \G^\ep$ and $\ell(\ep)$ replaced by $\widetilde{\G}^0, \widetilde{\G}^\ep$ and $\kappa^{-2}(\ep)$, respectively.
The verification of \textsc{\textbf{(A1)}}(i) will be shown in Proposition \ref{A11(MDP)} and  \textsc{\textbf{(A1)}}(ii) will be
presented in Proposition \ref{A12(MDP)}.
We first give a moment estimate for $Y^{\ep,h^\ep}$.

\beg{lem}\label{EsPerE(mdp)}
Suppose that  \textsc{\textbf{(H1)}} holds and let $\{h^\ep:\ep>0\}\subset\mathcal{A}_M$ for any $M\in(0,\infty)$.
Then for any $p\geq\th$, there exist two positive constants $C_{T,p,H,M}$ and $\ep_1$ such that
\begin{align*}
\sup_{\ep\in(0,\ep_1]}\E\bigg(\sup_{t\in[0,T]}|Y_t^{\ep,h^\ep}|^p\bigg)\leq C_{T,p,H,M}\bigg(1+\sup_{t\in[0,T]}|X_t^0|^p\bigg).
\end{align*}
\end{lem}
\beg{proof}
By \eqref{Meq-2}, we get
\begin{align}\label{1-PfEsPerE(mdp)}
|Y_t^{\ep,h^\ep}|^p\leq&\ff {3^{p-1}}{\ep^{pH}\kappa^p(\ep)}\left|\int_0^t\left(b(s, X^0_s+\ep^H\kappa(\ep)Y_s^{\ep, h^\ep}, \sL_{X_s^\ep})- b(s, X_s^0, \sL_{X_s^0})\right)\d s\right|^p\cr
&+3^{p-1}\left|\int_0^t\si(s,\sL_{X_s^\ep})\d(R_Hh^\ep)(s)\right|^p
+\ff{3^{p-1}}{\kappa^p(\ep)}\left|\int_0^t\si(s, \sL_{X_s^\ep})\d B_s^H\right|^p\cr
=:&I_1(t)+I_2(t)+I_3(t).
\end{align}
Using the H\"{o}lder inequality, \textsc{\textbf{(H1)}} and Lemma \ref{Dif(TE)} yields
\begin{align}\label{2-PfEsPerE(mdp)}
I_1(t)\leq&(6T)^{p-1}K^p(T)\int_0^t|Y_s^{\ep, h^\ep}|^p\d s+\ff {(6T)^{p-1}K^p(T)}{\ep^{pH}\kappa^p(\ep)}\int_0^t\W_\th^p( \sL_{X_s^\ep},\sL_{X_s^0})\d s\cr
\leq&(6T)^{p-1}K^p(T)\int_0^t|Y_s^{\ep, h^\ep}|^p\d s+\ff {(6T)^{p-1}K^p(T)}{\ep^{pH}\kappa^p(\ep)}\int_0^t\E|X_s^\ep-X_s^0|^p\d s\cr
\leq&(6T)^{p-1}K^p(T)\int_0^t|Y_s^{\ep, h^\ep}|^p\d s+\ff {6^{p-1}T^pK^p(T)C_{T,p,H}}{\kappa^p(\ep)}\bigg(1+\sup_{t\in[0,T]}|X_t^0|^p\bigg).
\end{align}
For the term $I_2(t)$, by \eqref{RFRH}, \textsc{\textbf{(H1)}}, Lemma \ref{Dif(TE)} and the Fubini theorem we obtain for any $p\geq\th$,
\begin{align}\label{3-PfEsPerE(mdp)}
I_2(t)
\leq&6^{p-1}\left|\int_0^t(\si(s,\sL_{X_s^\ep})-\si(s,\sL_{X_s^0}))\int_0^s\left(\ff s r\right)^{H-\ff 1 2}(s-r)^{H-\ff 3 2}(K_H^*h^\ep)(r)\d r\d s\right|^p\cr
&+6^{p-1}\left|\int_0^t\si(s,\sL_{X_s^0})\int_0^s\left(\ff s r\right)^{H-\ff 1 2}(s-r)^{H-\ff 3 2}(K_H^*h^\ep)(r)\d r\d s\right|^p\cr
\leq&C_{T,p,H}(1+\ep^{pH})\bigg(1+\sup_{t\in[0,T]}|X_t^0|^p\bigg)\cr
&\times\left(\int_0^tr^{\ff 1 2-H}|(K_H^*h^\ep)(r)|\left(\int_r^ts^{H-\ff 1 2}(s-r)^{H-\ff 3 2}\d s\right)\d r\right)^p\cr
\leq&C_{T,p,H}(1+\ep^{pH})\bigg(1+\sup_{t\in[0,T]}|X_t^0|^p\bigg)\|h^\ep\|_{\H}^p\cr
\leq&C_{T,p,H,M}(1+\ep^{pH})\bigg(1+\sup_{t\in[0,T]}|X_t^0|^p\bigg).
\end{align}
For the term $I_3(t)$, similar to \eqref{4-PfPrA12}, applying Lemma \ref{Dif(TE)} one sees that
\begin{align*}
\E\bigg(\sup_{t\in[0,T]}I_3(t)\bigg)\leq\ff{C_{T,p,H}(1+\ep^{pH})}{\kappa^p(\ep)}\bigg(1+\sup_{t\in[0,T]}|X_t^0|^p\bigg).
\end{align*}
Finally, plugging our previous inequalities \eqref{1-PfEsPerE(mdp)}-\eqref{3-PfEsPerE(mdp)} and resorting to the Gronwall inequality, we obtain
that there exist two positive constants $C_{T,p,H,M}$ and $\ep_1$ such that
\begin{align*}
\sup_{\ep\in(0,\ep_1]}\E\bigg(\sup_{t\in[0,T]}|Y_t^{\ep,h^\ep}|^p\bigg)\leq C_{T,p,H,M}\bigg(1+\sup_{t\in[0,T]}|X_t^0|^p\bigg),
\end{align*}
which completes the proof.
\end{proof}

\beg{prp}\label{A11(MDP)}
Suppose that  \textsc{\textbf{(H1)}} and  \textsc{\textbf{(H2)}} hold and let $\{h^\ep:\ep>0\}\subset\mathcal{A}_M$ for any $M\in(0,\infty)$.
Then for any $\delta>0$,
\begin{align*}
\lim_{\ep\ra 0}\P\left(\|Y_\cdot^{\ep, h^\ep}-\widetilde{\G}^0((R_Hh^\ep)(\cdot))\|_\infty>\delta\right)=0.
\end{align*}
\end{prp}
\beg{proof}
For each fixed $\ep>0$, using \eqref{Sk(mdp)}-\eqref{RateF(mdp)} with $h^\ep$ replacing $h$ and \eqref{Meq-2} yields
\begin{align*}
&Y_t^{\ep, h^\ep}-\widetilde{\G}^0((R_Hh^\ep)(t))=Y_t^{\ep, h^\ep}-\Xi_t^{h^\ep}\cr
=&\int_0^t\left[\ff 1{\ep^H\kappa(\ep)}\left(b(s, X^0_s+\ep^H\kappa(\ep)Y_s^{\ep, h^\ep}, \sL_{X_s^\ep})- b(s, X_s^0, \mathscr{L}_{X_s^0})\right)
-\nabla_{\Xi_s^{h^\ep}} b(s,\cdot,\sL_{X_s^0})(X_s^0)\right]\d s\cr
&+\int_0^t\left(\si(s,\sL_{X_s^\ep})-\sigma(s,\sL_{X_s^0})\right)\d(R_Hh^\ep)(s)+\ff 1{\kappa(\ep)}\int_0^t\si(s, \sL_{X_s^\ep})\d B_s^H,\ \ t\in[0,T].
\end{align*}
Then it follows that
\beg{align}\label{1-Pf(A11(MDP))}
&|Y_t^{\ep, h^\ep}-\Xi_t^{h^\ep}|^2\nonumber\\
\leq&3\bigg|\int_0^t\bigg[\ff 1{\ep^H\kappa(\ep)}\left(b(s, X^0_s+\ep^H\kappa(\ep)Y_s^{\ep, h^\ep}, \sL_{X_s^\ep})- b(s, X_s^0, \mathscr{L}_{X_s^0})\right)\cr
&\qquad\quad-\nabla_{\Xi_s^{h^\ep}} b(s,\cdot,\sL_{X_s^0})(X_s^0)\bigg]\d s\bigg|^2\nonumber\\
&+3\left|\int_0^t\left(\si(s,\sL_{X_s^\ep})-\sigma(s,\sL_{X_s^0})\right)\d(R_Hh^\ep)(s)\right|^2
+\ff 3{\kappa^2(\ep)}\left|\int_0^t\si(s, \sL_{X_s^\ep})\d B_s^H\right|^2\nonumber\\
=:&\sum_{i=1}^3J_i(t).
\end{align}
By \eqref{3-PfPrA12}, \eqref{4-PfPrA12} and Lemma \ref{Dif(TE)}, we get
\begin{align}\label{2-Pf(A11(MDP))}
J_2(t)\leq C_{T,H,M}\ep^{2H}\bigg(1+\sup_{t\in[0,T]}|X_t^0|^2\bigg)
\end{align}
and
\begin{align}\label{3-Pf(A11(MDP))}
\E\bigg(\sup_{t\in[0,T]}J_3(t)\bigg)\leq \ff{C_{T,H}(1+\ep^{2H})}{\kappa^2(\ep)}\bigg(1+\sup_{t\in[0,T]}|X_t^0|^2\bigg).
\end{align}
As for the term $J_1(t)$, applying the H\"{o}lder inequality,  \textsc{\textbf{(H1)}} and \textsc{\textbf{(H2)}} implies
\begin{align*}
&J_1(t)\cr
\leq&\ff{9T}{\ep^{2H}\kappa^2(\ep)}\int_0^t\left|\left(b(s, X^0_s+\ep^H\kappa(\ep)Y_s^{\ep, h^\ep}, \sL_{X_s^\ep})
-b(s, X^0_s+\ep^H\kappa(\ep)Y_s^{\ep, h^\ep}, \sL_{X_s^0})\right)\right|^2\d s\cr
&+9T\int_0^t\left|\ff 1{\ep^H\kappa(\ep)}\left(b(s, X^0_s+\ep^H\kappa(\ep)Y_s^{\ep, h^\ep}, \sL_{X_s^0})- b(s, X_s^0, \mathscr{L}_{X_s^0})\right)
-\nabla_{Y_s^{\ep,h^\ep}} b(s,\cdot,\sL_{X_s^0})(X_s^0)\right|^2\d s\cr
&+9T\int_0^t\left|\nabla_{Y_s^{\ep,h^\ep}} b(s,\cdot,\sL_{X_s^0})(X_s^0)
-\nabla_{\Xi_s^{h^\ep}} b(s,\cdot,\sL_{X_s^0})(X_s^0)\right|^2\d s\cr
\leq&\ff{9TK^2(T)}{\ep^{2H}\kappa^2(\ep)}\int_0^t\E\left|X_s^\ep-X_s^0\right|^2\d s+9T\widetilde{K}^2(T)\int_0^t|Y_s^{\ep, h^\ep}-\Xi_s^{h^\ep}|^2\d s\cr
&+9T\int_0^t\left|\int_0^1\left(\nabla_{Y_s^{\ep,h^\ep}} b(s,\cdot,\sL_{X_s^0})(Q_s^\ep(v))
-\nabla_{Y_s^{\ep,h^\ep}} b(s,\cdot,\sL_{X_s^0})(X_s^0)\right)\d v\right|^2\d s\cr
\leq&\ff{9TK^2(T)}{\ep^{2H}\kappa^2(\ep)}\int_0^t\E\left|X_s^\ep-X_s^0\right|^2\d s
+9T\widetilde{K}^2(T)\ep^{2H}\kappa^2(\ep)\int_0^t|Y_s^{\ep,h^\ep}|^4\d s\cr
&+9T\widetilde{K}^2(T)\int_0^t|Y_s^{\ep, h^\ep}-\Xi_s^{h^\ep}|^2\d s,
\end{align*}
where for any $v\in[0,1],Q_s^\ep(v)=X_s^0+v\ep^H\kappa(\ep)Y_s^{\ep, h^\ep}$.\\
Then, owing to Lemmas \ref{Dif(TE)} and \ref{EsPerE(mdp)}, we obtain
\begin{align}\label{4-Pf(A11(MDP))}
&\E\bigg(\sup_{s\in[0,t]}J_1(s)\bigg)\cr
\leq& \ff {C_{T,H}}{\kappa^2(\ep)}\bigg(1+\sup_{t\in[0,T]}|X_t^0|^2\bigg)+C_{T,H,M}\ep^{2H}\kappa^2(\ep)\bigg(1+\sup_{t\in[0,T]}|X_t^0|^4\bigg)\cr
&+9T\widetilde{K}^2(T)\int_0^t\sup_{r\in[0,s]}|Y_r^{\ep, h^\ep}-\Xi_r^{h^\ep}|^2\d s\cr
\leq&C_{T,H,M}\left(\ff 1{\kappa^2(\ep)}+\ep^{2H}\kappa^2(\ep)\right)\bigg(1+\sup_{t\in[0,T]}|X_t^0|^4\bigg)+C_T\int_0^t\sup_{r\in[0,s]}|Y_r^{\ep, h^\ep}-\Xi_r^{h^\ep}|^2\d s.
\end{align}
Gathering all the above estimates \eqref{2-Pf(A11(MDP))}-\eqref{4-Pf(A11(MDP))} into \eqref{1-Pf(A11(MDP))} and using the Gronwall inequality, we have thus obtained
\begin{align*}
\E\bigg(\sup_{t\in[0,T]}|Y_t^{\ep, h^\ep}-\Xi_t^{h^\ep}|^2\bigg)
\leq& C_{T,H,M}\left(\ep^{2H}+\ff{1+\ep^{2H}}{\kappa^2(\ep)}+\ff 1{\kappa^2(\ep)}+\ep^{2H}\kappa^2(\ep)\right)\bigg(1+\sup_{t\in[0,T]}|X_t^0|^4\bigg).
\end{align*}
Consequently, we arrive at
\begin{align*}
\lim_{\ep\ra0}\E\bigg(\sup_{t\in[0,T]}|Y_t^{\ep, h^\ep}-\Xi_t^{h^\ep}|^2\bigg)=0,
\end{align*}
which yields the desired assertion.
\end{proof}

\beg{prp}\label{A12(MDP)}
Suppose that  \textsc{\textbf{(H1)}} and  \textsc{\textbf{(H2)}} hold and let $\{h^n:n\in\mathbb{N}\}\subset\mathcal{S}_M$ for any $M\in(0,\infty)$
such that $h^n$ converges to element $h$ in $S_M$ as $n\ra\infty$.
Then
\begin{align*}
\lim_{n\ra\infty}\sup_{t\in[0,T]}|\widetilde{\G}^0(R_Hh^n)(t)-\widetilde{\G}^0(R_Hh)(t)|=0.
\end{align*}
\end{prp}
\beg{proof}
We pursue the same general strategy as in the same proof of Proposition \ref{PrA11}.
For every $n\geq1$, it follows from \eqref{RateF(mdp)} that $\widetilde{\G}^0(R_Hh)=\Xi^h$ and $\widetilde{\G}^0(R_Hh^n)=\Xi^{h^n}$,
where $\Xi^h$ (respectively $\Xi^{h^n}$) is the solution to \eqref{Sk(mdp)} (respectively  with $h^n$ replacing $h$).

We first show that $\{\Xi^{h^n}\}_{n\geq1}$ is relatively compact in $C([0,T]; \R^d)$.
Similar to Lemma \ref{Est(SkE)}, by \textsc{\textbf{(H2)}} it is easy to see that
\begin{align*}
\sup_{n\geq1}\sup_{t\in[0,T]}|\Xi_t^{h^n}|\leq C_{T,H,M},
\end{align*}
which implies that $\{\Xi^{h^n}\}_{n\geq1}$ is uniformly bounded in $C([0,T]; \R^d)$.
On the other hand, by \eqref{Sk(mdp)} and  \textsc{\textbf{(H2)}}, we have that for $0\leq s<t\leq T$,
\begin{align*}
|\Xi_t^{h^n}-\Xi_s^{h^n}|\leq&\left|\int_s^t\nabla_{\Xi_r^{h^n}} b(r,\cdot,\sL_{X_r^0})(X_r^0)\d r\right|+\left|\int_s^t\sigma(r,\sL_{X_r^0})\d(R_H h^n)(r)\right|\cr
\leq& C_{T,H,M}(t-s)+\left|\int_s^t\sigma(r,\sL_{X_r^0})\d(R_H h^n)(r)\right|.
\end{align*}
This, together with \eqref{3-PfPrA11}-\eqref{5-PfPrA11} and the fact that $\|h^n\|_\H\leq\sqrt{2M}$, leads to that $\{\Xi^{h^n}\}_{n\geq1}$ is equi-continuous in $C([0,T]; \R^d)$.
Therefore, we obtain relative compactness of $\{\Xi^{h^n}\}_{n\geq1}$ in $C([0,T]; \R^d)$ by applying the Arzel\`{a}-Ascoli theorem.
As a consequence, we can extract a  (not relabelled) subsequence such that $\Xi^{h^n}$ converges to some $\bar{\Xi}$ in $C([0,T]; \R^d)$.
Along the same lines as in the proof of Proposition \ref{PrA11}, we can obtain that $\bar{\Xi}=\Xi^h$, from which the result follows.
\end{proof}

\subsection{Proof of Theorem \ref{Th(clt)}}

\beg{proof}
For any $t\in[0,T]$, let $Z_t^\epsilon=\frac{X_t^\epsilon-X_t^0}{\epsilon^H}$.
By \eqref{Ap-DDsde} and \eqref{LimSDE}, it is easy to see that $Z_t^\epsilon$ satisfies
\beg{align*}
Z_t^\epsilon=\int_0^t\frac{1}{\epsilon^H}\Big(b(s, X_s^\epsilon, \mathscr{L}_{ X_s^\epsilon})-b(s, X_s^0, \mathscr{L}_{ X_s^0})\Big)\d s+\int_0^t\sigma(s,\mathscr{L}_{ X_s^\epsilon})\d B_s^H.
\end{align*}
This, together with \eqref{1-Th(clt)}, yields
\begin{align}\label{Pf0-Th(clt)}
&Z_t^\epsilon-Z_t\cr
=&\int_0^t\bigg(\frac{1}{\epsilon^H}\Big(b(s, X_s^\epsilon, \sL_{ X_s^\epsilon})-b(s, X_s^0, \sL_{X_s^0})\Big)-\nabla_{Z_s}b(s,\cdot ,\sL_{X_s^0})(X_s^0)\bigg)\d s\cr
&-\int_0^t\left(\E\langle D^Lb(s,u,\cdot)(\sL_{X_s^0})(X_s^0), Z_s\rangle\right)|_{u=X_s^0}\d s\cr
&+\int_0^t\left(\sigma(s,\sL_{ X_s^\epsilon})-\sigma(s,\sL_{ X_s^0})\right)\d B_s^H\cr
= &\int_0^t\bigg[\frac{1}{\epsilon^H}\Big(b(s, X_s^\epsilon, \sL_{ X_s^\epsilon})-b(s, X_s^0,\sL_{ X_s^0})\Big)\cr
&\ \ \ \ \ \ -\nabla_{Z_s^\epsilon}b(s,\cdot ,\sL_{X_s^\epsilon})(X_s^0)
-\left(\E\langle D^Lb(s,u,\cdot)(\sL_{X_s^0})(X_s^0), Z_s^\ep\rangle\right)|_{u=X_s^0}\bigg]\d s\cr
& +\int_0^t\left(\nabla_{Z_s^\epsilon}b(s,\cdot ,\sL_{X_s^\epsilon})(X_s^0)-\nabla_{Z_s}b(s,\cdot,\sL_{X_s^0})(X_s^0)\right)\d s\cr
& +\int_0^t\left(\E\langle D^Lb(s,u,\cdot)(\sL_{X_s^0})(X_s^0), Z_s^\ep-Z_s\rangle\right)|_{u=X_s^0}\d s\cr
& +\int_0^t\left(\sigma(s,\sL_{ X_s^\epsilon})-\sigma(s,\sL_{ X_s^0})\right)\d B_s^H\cr
=:&\sum_{i=1}^{4}I_i(t).
\end{align}
Applying Lemma \ref{FoLD} and \textsc{\textbf{(H3)}}(ii), we get
\begin{align*}
&|I_1(t)|\cr
=&\bigg|\int_0^t\bigg[\int_0^1\left(\frac{1}{\ep^H}\frac{\d}{\d\tilde{\th}}b(s,X_s^{\ep,0}(\tilde{\th}), \sL_{ X_s^\ep})
-\nabla_{Z_s^\ep}b(s,\cdot ,\sL_{X_s^\ep})(X_s^0)\right)\d\tilde{\th}\\
&\ \ \ \ \  +\int_0^1\left(\frac{1}{\ep^H}\frac{\d}{\d\tilde{\th}}b(s,X_s^0,\sL_{X_s^{\ep,0}(\tilde{\th})})
-\left(\E\langle D^Lb(s,u,\cdot)(\sL_{X_s^0})(X_s^0), Z_s^\ep\rangle\right)|_{u=X_s^0}\right)
\d\tilde{\th}\bigg]\d s\bigg|\\
=&\bigg|\int_0^t\bigg[\int_0^1\left(\nabla_{Z_s^\ep} b(s,\cdot, \sL_{ X_s^\ep})(X_s^{\ep,0}(\tilde{\th}))
-\nabla_{Z_s^\ep}b(s,\cdot ,\sL_{X_s^\ep})(X_s^0)\right)\d\tilde{\th}\\
&\ \ \ \ \ +\int_0^1\left(\E\langle D^Lb(s,v,\cdot)(\sL_{X_s^{\ep,0}(\tilde{\th})})(X_s^{\ep,0}(\tilde{\th}))-D^Lb(s,u,\cdot)(\sL_{X_s^0})(X_s^0), Z_s^\ep\rangle\right)|_{v=X_s^0,u=X_s^0}\d\tilde{\th}\bigg]\d s\bigg|\\
\leq&\int_0^t|Z_s^\ep|\int_0^1\bar{K}(s)\left(|X_s^{\ep,0}(\tilde{\th})-X_s^0|+\E|X_s^{\ep,0}(\tilde{\th})-X_s^0|+\W_\th(\sL_{X_s^{\ep,0}(\tilde{\th})},\sL_{X_s^0})\right)\d\tilde{\th}\d s\cr
\leq&\bar{K}(t)\int_0^t|Z_s^\ep|\left(|X_s^\ep-X_s^0|+\E|X_s^\ep-X_s^0|+\left(\E|X_s^\ep-X_s^0|^\th\right)^{\ff 1 \th}\right)\d s,
\end{align*}
where for any $\tilde{\th}\in[0,1],X_s^{\ep,0}(\tilde{\th})=X_s^0+\tilde{\th}(X_s^\ep-X_s^0)$.\\
Consequently, by Lemma \ref{Dif(TE)} we have for any $p\geq\th$,
\begin{align}\label{Pf1-Th(clt)}
\E\bigg(\sup_{0\leq t\leq T}|I_1(t)|^p\bigg)\le C_{T,p,H}\ep^{pH}\left(1+\sup_{t\in[0,T]}|X_t^0|^{2p}\right).
\end{align}
For the terms $I_i(t),i=2,3$, from \textsc{\textbf{(H3)}} again we have
\begin{align*}
|I_2(t)|
\leq&\int_0^t\Big(\left|\nabla_{Z_s^\epsilon}b(s,\cdot ,\mathscr{L}_{X_s^\epsilon})(X_s^0)-\nabla_{Z_s^\epsilon}b(s,\cdot ,\mathscr{L}_{X_s^0})(X_s^0)\right|\cr
&\ \ \ \ \ \ \ +\left|\nabla_{Z_s^\epsilon}b(s,\cdot ,\mathscr{L}_{X_s^0})(X_s^0)-\nabla_{Z_s}b(s,\cdot ,\mathscr{L}_{X_s^0})(X_s^0)\right|\Big)\d s\cr
\leq&\bar{K}(t)\int_0^t\left(|Z_s^\epsilon|\W_\th(\sL_{X_s^\epsilon},\mathscr{L}_{X_s^0})+|Z_s^\epsilon-Z_s|\right)\d s
\end{align*}
and
\begin{align*}
|I_3(t)|\leq\bar{K}(t)\int_0^t\E|Z_s^\epsilon-Z_s|\d s.
\end{align*}
Then, Lemma \ref{Dif(TE)} implies
\begin{align}\label{Pf2-Th(clt)}
&\E\bigg(\sup_{0\leq t\leq T}|I_2(t)+I_3(t)|^p\bigg)\cr
\leq& C_{T,p,H}\left[\ep^{pH}\bigg(1+\sup_{t\in[0,T]}|X_t^0|^{2p}\bigg)+\int_0^T\E\bigg(\sup_{r\in[0,s]}|Z_r^\epsilon-Z_r|^p\bigg)\d s\right].
\end{align}
As for the term $I_4(t)$, by Lemma \ref{MomEs} and \textsc{\textbf{(H3)}}(i), we obtain for any $p>1/H$ and $p\geq\th$,
\begin{align}\label{Pf3-Th(clt)}
\E\bigg(\sup_{0\leq t\leq T}|I_4(t)|^p\bigg)\leq& C_{T,p,H}\int_0^T\left\|\si(s,\sL_{X_s^\ep})-\si(s,\sL_{X_s^0})\right\|^p\d s\cr
\leq&C_{T,p,H}\int_0^T\bar{K}^p(s)\W_\th^p(\sL_{X_s^\ep},\sL_{X_s^0})\d s\cr
\leq&C_{T,p,H}\int_0^T\E|X_s^\ep-X_s^0|^p\d s\cr
\leq&C_{T,p,H}\ep^{pH}\left(1+\sup_{t\in[0,T]}|X_t^0|^p\right),
\end{align}
where the last inequality is due to Lemma \ref{Dif(TE)}.\\
Plugging \eqref{Pf1-Th(clt)}-\eqref{Pf3-Th(clt)} into \eqref{Pf0-Th(clt)}, we end up with
\begin{align*}
\E\bigg(\sup_{0\leq t\leq T}|Z_t^\epsilon-Z_t|^p\bigg)\leq
C_{T,p,H}\left[\ep^{pH}\bigg(1+\sup_{t\in[0,T]}|X_t^0|^{2p}\bigg)+\int_0^T\E\bigg(\sup_{r\in[0,s]}|Z_r^\epsilon-Z_r|^p\bigg)\d s\right].
\end{align*}
Therefore, a simple application of the Gronwall inequality yields the desired result.
\end{proof}

\section{Appendix: proofs of the auxiliary results}

\subsection{Proof of Lemma \ref{Le(Meas)}}

For any $\mu\in C([0,T]; \sP_2(\R^d))$, we let
\begin{align*}
b^\mu(t,x)=b(t,x,\mu_t),\ \  \si^\mu(t)=\si(t,\mu_t), \ \ t\in[0,T], \ x\in\R^d.
\end{align*}
It is readily checked that the functions $b^\mu(t,x)$ and $\si^\mu(t)$ also satisfy \textsc{\textbf{(H1)}}.
Then, by \cite[Theorem 3.1]{FHSY}, it is known that equation \eqref{FrozeEQ} has a unique solution $\tilde{X}\in\cS^p([0,T])$ with $p>1/H$.
As a consequence, there is a measurable map $\G_\mu: C([0,T]; \R^d)\rightarrow C([0,T]; \R^d)$ such that
\begin{align*}
\tilde{X}_\cdot=\mathcal{G}_\mu(B_\cdot^H).
\end{align*}

Next, we focus on proving the other assertion.

For any $h\in\mathcal{A}_M$,  by \eqref{IRFor} and \eqref{RHop} we have
\begin{align*}
\tilde{B}_\cdot^H:=B_\cdot^H+(R_Hh)(\cdot)=\int_0^\cdot K_H(\cdot,s)(\d W_s+(K_H^*h)(s)\d s).
\end{align*}
Set
\begin{align*}
\vartheta_T:=\exp\left[-\int_0^T\langle(K_H^*h)(s),\d W_s\rangle-\ff 1 2\int_0^T|(K_H^*h)(s)|^2\d s\right].
\end{align*}
Taking into account of the isometry of $K_H^*$  between $\H$ and $L^2([0,T],\R^d)$, we obtain
\begin{align*}
\E\exp\left[\ff 1 2\int_0^T|(K_H^*h)(s)|^2\d s\right]=\E\exp\left[\ff 1 2\|h\|^2_\H\right]\leq\e^M.
\end{align*}
Then, according to the Girsanov theorem for the fractional Brownian motion (see, e.g., \cite[Theorem 4.9]{Decreusefond&Ustunel98a} or \cite[Theorem 2]{Nualart&Ouknine02b}),
$\{\tilde{B}_t^H\}_{t\in[0,T]}$ is a $d$-dimensional fractional Brownian motion under the probability $\vartheta_T\P$.
So, we deduce that
\begin{align*}
\tilde{X}_\cdot^h:=\mathcal{G}_\mu\left(B_\cdot^H+(R_Hh)(\cdot)\right)=\mathcal{G}_\mu(\tilde{B}_\cdot^H)
\end{align*}
is the unique solution of the equation \eqref{FrozeEQ} on $(\Omega,\sF,\{\sF_t\}_{t\in[0,T]},\vartheta_T\P)$.
Since $\P$ and $\vartheta_T\P$ are two equivalent probability measures, on can see that $\tilde{X}_\cdot^h$ is a strong solution of  the equation \eqref{1-Le(Meas)} on $(\Omega,\sF,\{\sF_t\}_{t\in[0,T]},\P)$.
Along the same lines as above, we can also derive the uniqueness of the solution $\tilde{X}_\cdot^h$.
The proof is completed.\qed

\subsection{Proof of Lemma \ref{MomEs}}

Before proving Lemma \ref{MomEs}, we first introduce the Hardy-Littlewood inequality (see, e.g., \cite[Theroem 1, Page 119]{Stein70}).

\beg{lem}\label{HLI}
Let $1<\tilde{p}<\tilde{q}<\infty$ and $\frac{1}{\tilde{q}}=\frac{1}{\tilde{p}}-\alpha$.
Suppose that $f:\R_+\rightarrow\R$ is in $L^{\tilde{p}}(0,\infty)$, then $I^\alpha_{0+}f(x)$ converges absolutely for almost every $x$.
Furthermore, there exists some positive constant $C_{\tilde{p},\tilde{q}}$ such that
\begin{align*}
\|I^\alpha_{0+}f\|_{L^{\tilde{q}}(0,\infty)}\leq C_{\tilde{p},\tilde{q}}\|f\|_{L^{\tilde{p}}(0,\infty)},
\end{align*}
where $I^\alpha_{0+}f$ is the left-sided  fractional Riemann-Liouville integral of $f$ of order $\alpha$ defined as
\begin{align}\label{1-HLI}
&I_{0+}^\alpha f(x)=\frac{1}{\Gamma(\alpha)}\int_0^x\frac{f(y)}{(x-y)^{1-\alpha}}\d y.
\end{align}
\end{lem}

\emph{Proof of Lemma \ref{MomEs}.}
Owing to \textsc{\textbf{(H1)}} and $\mu\in C([0,T]; \sP_p(\R^d))$ with $p\geq\th$, one can verify that $\si(\cdot,\mu_\cdot)\in L^p([0,T],\R^d\otimes\R^d)$.
Then by the H\"{o}lder inequality and  $p>1/H$, there holds $\si(\cdot,\mu_\cdot)\in L^{\ff 1 H}([0,T],\R^d\otimes\R^d)$,
which means that for any $t\in[0,T], \int_0^t\sigma (s,\mu_s)\d B_s^H $ is well-defined because of $ L^{\ff 1 H}([0,T],\R^d\otimes\R^d)\subset\H$.

Next, we intend to show \eqref{MomEs-1}.

Since $pH>1$, we choose $\la$ such that $1-H<\la<1-1/p$ and put $C_\la:=\int_s^t(t-r)^{-\la}(r-s)^{\la-1}\d r$.
Applying the stochastic Fubini theorem and the H\"{o}lder inequality, we have
\begin{align}\label{Pf1-MomEs}
&\E\left(\sup\limits_{t\in[0,T]}\left|\int_0^t\si(s,\mu_s)\d B_s^H\right|^p\right)\cr
&=C_\la^{-p}\E\left(\sup\limits_{t\in[0,T]}\left|\int_0^t\left(\int_s^t(t-r)^{-\la}(r-s)^{\la-1}\d r\right)\si(s,\mu_s)\d B_s^H\right|^p\right)\cr
&=C_\la^{-p}\E\left(\sup\limits_{t\in[0,T]}\left|\int_0^t (t-r)^{-\la}\left(\int_0^r(r-s)^{\la-1}\si(s,\mu_s)\d B_s^H\right)\d r\right|^p\right)\cr
&\leq\ff {C_\la^{-p}(p-1)^{p-1}}{(p-1-\la p)^{p-1}}\E\left(\sup\limits_{t\in[0,T]}t^{p-1-\la p}
\int_0^t \left|\int_0^r(r-s)^{\la-1}\si(s,\mu_s)\d B_s^H\right|^p\d r\right)\cr
&\leq\ff {C_\la^{-p}(p-1)^{p-1}}{(p-1-\la p)^{p-1}}T^{p-1-\la p}\int_0^T\E\left|\int_0^r(r-s)^{\la-1}\si(s,\mu_s)\d B_s^H\right|^p\d r.
\end{align}
Here we have used the condition $\la<1-1/p$ in the first inequality.\\
Observe that for each $r\in[0,T], \int_0^r(r-s)^{\la-1}\si(s,\sL_{X^n_s})\d B_s^H$ is a centered Gaussian random variable.
Consequently, with the help of the Kahane-Khintchine formula, we derive that there is a constant $C_p>0$ such that
\begin{align*}
&\E\left|\int_0^r(r-s)^{\la-1}\si(s,\mu_s)\d B_s^H\right|^p\cr
&\leq C_p\left(\E\left|\int_0^r(r-s)^{\la-1}\si(s,\mu_s)\d B_s^H\right|^2\right)^{\ff p 2}\cr
&\leq C_{p,H}\left(\int_0^r\int_0^r(r-u)^{\la-1}\|\si(u,\mu_u)\|(r-v)^{\la-1}\|\si(v,\mu_v)\|\cdot|u-v|^{2H-2}\d u\d v\right)^{\ff p 2}\cr
&\leq C_{p,H}\left(\int_0^r(r-s)^{\ff {\la-1} H} \|\si(s,\mu_s)\|^{\ff 1 H}\d s\right)^{pH}.
\end{align*}
Here, we have adopted the argument in \cite[Theroem 1.1, Page 201]{JMV01} in the last inequality. \\
Then, plugging this into \eqref{Pf1-MomEs} and taking into account of the condition $1-H<\la$ and Lemma \ref{HLI} with $\tilde{q}=pH$ and $\al=1-\ff{1-\la}H$ (implying $\tilde{p}=\ff{pH}{p(\la+H-1)+1}$),  we conclude that
\begin{align*}
&\E\left(\sup\limits_{t\in[0,T]}\left|\int_0^t\si(s,\mu_s)\d B_s^H\right|^p\right)\cr
&\leq C_{\la,p,H}T^{p-1-\la p}\int_0^T\left(\int_0^r(r-s)^{\ff {\la-1} H} \|\si(s,\mu_s)\|^{\ff 1 H}\d s\right)^{pH}\d r\cr
&\leq C_{\la,p,H}T^{p-1-\la p}\left(\int_0^T \|\si(r,\mu_s)\|^{\ff p {p(\la+H-1)+1}}\d r\right)^{p(\la+H-1)+1}\cr
&\leq C_{\la,p,H}T^{pH-1}\int_0^T\|\si(s,\mu_s)\|^p\d s,
\end{align*}
where the last inequality is due to the H\"{o}lder inequality.\\
Note that by taking proper $\la$, the constant $C_{\la,p,H}$ above may depend only on $p$ and $H$, which completes the proof.
\qed

\subsection{Proof of Proposition \ref{Suf1(LDP)}}

The proof of Proposition \ref{Suf1(LDP)} follows the method in \cite[Theorem 4.4]{BD00} in general.
It relies on the following preliminary result concerning a variational representation for random functional,
which is a slight change of \cite[Theorem 3.2]{Zhang09}.

\beg{lem}\label{VaRe}
Let $f$ be a bounded Borel measurable function on $\Omega$. Then there holds
\begin{align*}
-\log\E(\e^{-f})=\inf_{h\in\A}\E\left(f(\cdot+R_Hh)+\ff 1 2\|h\|_\H^2\right).
\end{align*}
\end{lem}

\emph{Proof of Proposition \ref{Suf1(LDP)}.}
In order to prove the proposition, we need to show that
 \eqref{De(lp)1} and  \eqref{De(lp)2} holds for all real-valued, bounded and continuous function $\varrho$ on $\mathscr{E}$,
and $I$ given in \eqref{1-Suf1(LDP)} is a rate function.
Notice first that the fact that $I$ is a rate function is readily checked via  \textsc{\textbf{(A0)}}(ii),
and it is thus omitted here for the sake of conciseness.
Below we shall focus on handling \eqref{De(lp)1} and  \eqref{De(lp)2}.

By Lemma \ref{VaRe} with $f(\cdot)$ replaced by $\ff {\varrho\circ\G^\ep(\ep^H\cdot)}{\ell(\ep)}$, we have
\begin{align}\label{1-PfSuf1(LDP)}
-\ell(\ep)\log\E\left[\exp\left(-\ff {\varrho(\mathbb{X}^\ep)}{\ell(\ep)} \right)\right]
=&-\ell(\ep)\log\E\left[\exp\left(-\ff {\varrho\circ\G^\ep(\ep^H B_\cdot^H)}{\ell(\ep)} \right)\right]\cr
=&\inf_{h\in\A}\E\left(\varrho\circ\G^\ep(\ep^H(B^H+R_Hh))+\ff 1 2\ell(\ep)\|h\|_\H^2 \right)\cr
=&\inf_{h\in\A}\E\left(\varrho\circ\G^\ep(\ep^H(B^H+R_Hh/{\ell^{\ff 1 2}(\ep)}))+\ff 1 2\|h\|_\H^2 \right).
\end{align}
The rest of the proof will be divided into two steps.

\textsl{Step 1. The upper bound}.
Without lost of generality, we assume that $\inf_{x\in\mathscr{E}}\{\varrho(x)+I(x)\}<\infty$.
Let $\delta>0$ be fixed. Then there exists $x_0\in\mathscr{E}$ such that
\begin{align}\label{2-PfSuf1(LDP)}
\varrho(x_0)+I(x_0)\leq\inf_{x\in\mathscr{E}}\{\varrho(x)+I(x)\}+\ff \delta 2.
\end{align}
By \eqref{1-Suf1(LDP)}, we choose $\hbar\in\H$ such that $\G^0(R_H\hbar)=x_0$ and
\begin{align*}
\ff 1 2 \|\hbar\|_{\H}^2\leq  I(x_0)+\ff \delta 2.
\end{align*}
Then, combining this  with \eqref{1-PfSuf1(LDP)} yields
\begin{align*}
-\ell(\ep)\log\E\left[\exp\left(-\ff {\varrho(\mathbb{X}^\ep)}{\ell(\ep)} \right)\right]
\leq \E\left(\varrho\circ\G^\ep(\ep^H(B^H+R_H\hbar/\ell^{\ff 1 2}(\ep)))\right)+I(x_0)+\ff \delta 2.
\end{align*}
Since $\varrho$ is bounded and continuous, taking limit as $\ep\ra0$ and using  \textsc{\textbf{(A0)}}(i) imply
\begin{align*}
&\limsup_{\ep\ra0}-\ell(\ep)\log\E\left[\exp\left(-\ff {\varrho(\mathbb{X}^\ep)}{\ell(\ep)} \right)\right]\cr
\leq& \varrho\circ\G^0(R_H\hbar))+I(x_0)+\ff \delta 2\cr
=& \varrho(x_0)+I(x_0)+\ff \delta 2\cr
\leq&\inf_{x\in\mathscr{E}}\{\varrho(x)+I(x)\}+\delta,
\end{align*}
where the last inequality follows from \eqref{2-PfSuf1(LDP)}.
Since $\delta>0$ is arbitrary, we complete the proof of the upper bound.

\textsl{Step 2. The lower bound}.
Fix $\delta>0$.
According to \eqref{1-PfSuf1(LDP)}, for every $\ep>0$ there exists $h^\ep\in\A$ such that
\begin{align}\label{3-PfSuf1(LDP)}
-\ell(\ep)\log\E\left[\exp\left(-\ff {\varrho(\mathbb{X}^\ep)}{\ell(\ep)} \right)\right]
\geq\E\left(\varrho\circ\G^\ep(\ep^H(B^H+R_Hh^\ep/{\ell^{\ff 1 2}(\ep)}))+\ff 1 2\|h^\ep\|_\H^2 \right)-\de,
\end{align}
which also implies
\begin{align}\label{4-PfSuf1(LDP)}
\sup_{\ep>0}\E\left(\ff 1 2\|h^\ep\|_\H^2 \right)\leq2\|\varrho\|_\infty+\de.
\end{align}
Now, for each finite number $N$, we define stopping times
\begin{align*}
\tau_N^\ep:=\inf\left\{t\in[0,T]: \ff 1 2\|h^\ep\mathrm{I}_{[0,t]}\|_\H^2\geq N\right\}\wedge T,
\end{align*}
and the processes $h^{\ep,N}(t):=h^\ep(t)\mathrm{I}_{[0,\tau_N^\ep]}(t)$.
One can see that $h^{\ep,N}\in\A$ and
\begin{align}\label{5-PfSuf1(LDP)}
\P\left(h^\ep\neq h^{\ep,N}\right)\leq\P\left(\ff 1 2\|h^\ep\|_\H^2\geq N\right)\leq\ff{2\|\varrho\|_\infty+\de}{N},
\end{align}
where the last inequality is due to the Markov inequality and \eqref{4-PfSuf1(LDP)}.\\
Moreover, observe that we have
\begin{align}\label{6-PfSuf1(LDP)}
&\varrho\circ\G^\ep(\ep^H(B^H+R_Hh^\ep/{\ell^{\ff 1 2}(\ep)}))\cr
=&\varrho\circ\G^\ep(\ep^H(B^H+R_Hh^{\ep,N}/{\ell^{\ff 1 2}(\ep)}))\cr
&+\left[\varrho\circ\G^\ep(\ep^H(B^H+R_Hh^\ep/{\ell^{\ff 1 2}(\ep)}))-\varrho\circ\G^\ep(\ep^H(B^H+R_Hh^{\ep,N}/{\ell^{\ff 1 2}(\ep)}))\right]\mathrm{I}_{\{h^\ep\neq h^{\ep,N}\}}\cr
\geq&\varrho\circ\G^\ep(\ep^H(B^H+R_Hh^{\ep,N}/{\ell^{\ff 1 2}(\ep)}))-2\|\varrho\|_\infty\mathrm{I}_{\{h^\ep\neq h^{\ep,N}\}},
\end{align}
and by \eqref{Isom} we get
\begin{align}\label{7-PfSuf1(LDP)}
\|h^\ep\|_\H^2=\|K_H^*h^\ep\|_{L^2}\geq\|K_H^*h^{\ep,N}\|_{L^2}=\|h^{\ep,N}\|_\H^2.
\end{align}
Plugging \eqref{5-PfSuf1(LDP)}-\eqref{7-PfSuf1(LDP)} into \eqref{3-PfSuf1(LDP)} yields
\begin{align}\label{8-PfSuf1(LDP)}
&-\ell(\ep)\log\E\left[\exp\left(-\ff {\varrho(\mathbb{X}^\ep)}{\ell(\ep)} \right)\right]\cr
\geq&\E\left(\varrho\circ\G^\ep(\ep^H(B^H+R_Hh^{\ep,N}/{\ell^{\ff 1 2}(\ep)}))+\ff 1 2\|h^{\ep,N}\|_\H^2 \right)-\ff{2\|\varrho\|_\infty(2\|\varrho\|_\infty+\de)}{N}-\de.
\end{align}
Since $N$ and $\de$ are arbitrary in \eqref{8-PfSuf1(LDP)}, in proving the lower bound it suffices to show that
\begin{align}\label{9-PfSuf1(LDP)}
\liminf_{\ep\ra0}\E\left(\varrho\circ\G^\ep(\ep^H(B^H+R_Hh^{\ep,N}/{\ell^{\ff 1 2}(\ep)}))+\ff 1 2\|h^{\ep,N}\|_\H^2 \right)\geq\inf_{x\in\mathscr{E}}\{\varrho(x)+I(x)\}.
\end{align}
By the definition of $h^{\ep,N}$, one has
\begin{align*}
\sup_{\ep>0}\ff 1 2\|h^{\ep,N}\|_\H^2\leq N,\  a.s.,
\end{align*}
which allows us to extract a (not relabelled) subsequence such that $h^{\ep,N}$ converges to $h$ in distribution as $S_N$-valued random elements.
Then, by  \textsc{\textbf{(A0)}}(i) and the Fubini theorem we obtain
\begin{align*}
&\liminf_{\ep\ra0}\E\left(\varrho\circ\G^\ep(\ep^H(B^H+R_Hh^{\ep,N}/{\ell^{\ff 1 2}(\ep)}))+\ff 1 2\|h^{\ep,N}\|_\H^2 \right)\cr
\geq& \E\left(\varrho(\G^0(R_Hh))+\ff 1 2\|h\|_\H^2 \right)\cr
\geq& \inf_{\{(x,h)\in\mathscr{E}\times\H:x=\G^0(R_Hh)\}}\left\{\varrho(x)+\ff 1 2\|h\|_\H^2\right\}\cr
\geq& \inf_{x\in\mathscr{E}}\{\varrho(x)+I(x)\},
\end{align*}
which completes the proof of the lower bound.
The proof is therefore finished.
\qed

\textbf{Acknowledgement}

X. Fan is partially supported by the Natural Science Foundation of Anhui Province (No. 2008085MA10) and the National Natural Science Foundation of China (No. 11871076, 12071003).


\begin{thebibliography}{17}
{\small

\setlength{\baselineskip}{0.14in}
\parskip=0pt



\bibitem{Alos&Mazet&Nualart01a} E. Al\`{o}s, O. Mazet and D. Nualart, Stochastic calculus with respect to Gaussian processes, \textit{Ann. Probab.} {\bf 29} (2001), 766--801.


\bibitem{Biagini&Hu08a} F. Biagini, Y. Hu, B. $\emptyset$ksendal and  T. Zhang, Stochastic Calculus for Fractional Brownian Motion and Applications, Springer-Verlag, London, 2008.


\bibitem{BT97} M. Bossy and D. Talay, A stochastic particle method for the McKean-Vlasov and the Burgers equation, \textit{Math. Comput.} {\bf 66} (1997), 157--192.

\bibitem{BDS22} S. Bourguin, T. Dang, and K. Spiliopoulos, Moderate deviation principle for multiscale systems driven by fractional Brownian motion,
\textit{arXiv:2206.06794}.

\bibitem{BGJ17} Z. Brzeniak,  B. Goldys and T. Jegaraj, Large deviations and transitions between equilibria for stochastic Landau-Lifshitz-Gilbert equation, \textit{Arch. Ration. Mech. Anal.} {\bf 226} (2017), 497--558.


\bibitem{Buckdahn&L&Peng&ainer17a} R. Buckdahn, J. Li, S. Peng and C. Rainer,  Mean-field stochastic differential equations and
associated PDEs, \textit{Ann. Probab.} {\bf 2} (2017), 824--878.


\bibitem{BD00} A. Budhiraja and P. Dupuis, A variational representation for positive functionals of infinite dimensional Brownian motion, \textit{Probab. Math. Statist.-Wroclaw University} {\bf 20} (2000), 39--61.


\bibitem{BD19} A. Budhiraja and P. Dupuis, Analysis and Approximation of Rare Events: Representations and Weak Convergence Methods, Springer, 2019.


\bibitem{BDG16}  A. Budhiraja, P. Dupuis and A. Ganguly,  Moderate deviation principles for stochastic differential equations with jumps, \textit{ Ann. Probab.} {\bf 44} (2016), 1723--1775.


\bibitem{BDM08} A. Budhiraja, P. Dupuis and V. Maroulas, Large deviations for infinite dimensional stochastic dynamical systems, \textit{Ann. Probab.} {\bf 36} (2008), 1390--1420.

\bibitem{BDM11} A. Budhiraja, P. Dupuis and V. Maroulas, Variational representations for continuous time processes, \textit{Ann. Inst. H. Poincar\'{e} Probab. Statist.} {\bf 47} (2011), 725--747.

\bibitem{BS20} A. Budhiraja, and X. Song, Large deviation principles for stochastic dynamical systems with a fractional Brownian noise,
\textit{arXiv:2006.07683}.




\bibitem{Cardaliaguet13} P. Cardaliaguet, Notes on mean field games, P.-L. Lions lectures at Coll\`{e}ge de France,
https://www.ceremade.dauphine.fr/cardaliaguet/MFG20130420.pdf, 2013.

\bibitem{CD13} R. Carmona and F. Delarue, Probabilistic analysis of mean-field games, \textit{SIAM J. Control Optim.} {\bf 51} (2013), 2705--2734.

\bibitem{Crisan&McMurray18a} D. Crisan and E. McMurray, Smoothing properties of McKean-Vlasov SDEs,
\textit{Probab. Theory Related Fields} {\bf 171} (2018), 97--148.




\bibitem{Decreusefond&Ustunel98a}L. Decreusefond and A. S. \"{U}st\"{u}nel, Stochastic analysis of the fractional Brownian motion, \textit{Potential Anal.} {\bf 10} (1998), 177--214.

\bibitem{DST19} G. Dos Reis, W. Salkeld and J. Tugaut, Freidlin-Wentzell LDP in path space for McKean-Vlasov equations and the functional iterated logarithm law, \textit{Ann. Appl. Probab.} {\bf 29} (2019), 1487--1540.

\bibitem{DWZZ20} Z. Dong, J. Wu, R. Zhang and T. Zhang, Large deviation principles for first-order scalar conservation laws with stochastic forcing, \textit{Ann. Appl. Probab.} {\bf 30} (2020), 324--367.

\bibitem{DXZZ17} Z. Dong, J. Xiong, J. Zhai and T. Zhang, A moderate deviation principle for 2-D stochastic Navier-Stokes equations driven by multiplicative l\'{e}vy noises, \textit{J. Funct. Anal.} {\bf 272} (2017), 227--254.

\bibitem{DE11} P. Dupuis and R. Ellis, A Weak Convergence Approach to the Theory of Large Deviations, John Wiley  Sons, 2011.





\bibitem{FHSY} X. Fan, X. Huang, Y. Suo and C. Yuan, Distribution dependent SDEs driven by fractional Brownian motions, \textit{Stochastic Process. Appl.}
{\bf 151} (2022), 23--67.

\bibitem{FW84} M. I. Freidlin and A. D. Wentzell, Random perturbations of dynamical systems, \textit{Grundlehren der Mathematischen Wissenschaften [Fundamental Principles of Mathematical Sciences]}, Volume 260, Springer, 1984.




\bibitem{HLL21} W. Hong, S. Li and W. Liu, Large deviation principle for McKean-Vlasov quasilinear stochastic evolution equations,
\textit{Appl. Math. Optim.} {\bf 84} (2021), 1119--1147.

\bibitem{HW} X. Huang and F.-Y. Wang, Distribution dependent SDEs with singular coefficients, \textit{Stochastic Process. Appl.} {\bf 129} (2019), 4747--4770.

\bibitem{HW21b} X. Huang and F.-Y. Wang, Mckean-Vlasov SDEs with drifts discontinuous under wasserstein distance, \textit{Discrete Contin. Dyn. Syst.} {\bf 41} (2021), 1667--1679.




\bibitem{JW17} P. E. Jabin and Z. Wang, Mean field limit for stochastic particle systems,  In \textit{Active Particles},
Volume 1, pages 379--402. Springer, 2017.

\bibitem{JY21} H. Jiang and Q. Yang, Sample path large deviations for the multiplicative Poisson shot noise process with compensation, \textit{Stochastics} {\bf 93} (2021), 447--477.




\bibitem{LL07} J. Lasry and P. Lions, Mean field games, \textit{Jpn. J. Math.} {\bf 2} (2007), 229--260.

\bibitem{LWbook} W. Liu and M. R\"{o}ckner, Stochastic Partial Differential Equations: An Introduction, Universitext, Springer, 2015.

\bibitem{LSZZ22} W. Liu, Y. Song, J. Zhai and T. Zhang, Large and moderate deviation principles for McKean-Vlasov SDEs with jumps, \textit{Potential Anal.} (2022), 1--50.




\bibitem{MSZ21} A. Matoussi, W. Sabbagh and T. Zhang, Large deviation principle of obstacle problems for quasilinear stochastic PDEs, \textit{Appl. Math. Optim.} {\bf 83} (2021), 849--879.

\bibitem{McKean66} H. P. McKean, A class of Markov processes associated with nonlinear parabolic equations, \textit{Proc. Natl. Acad. Sci. USA} {\bf 56} (1966), 1907--1911.

\bibitem{JMV01} J. M\'{e}min, Y. Mishura and E. Valkeila, Inequalities for the moments of Wiener integrals with respect to a fractional Brownian motion,
\textit{Statist. Probab. Lett.} {\bf 51} (2001), 197--206.



\bibitem{ND} D. Nualart, The Malliavin Calculus and Related Topics, Second edition, Springer-Verlag, Berlin, 2006.

\bibitem{Nualart&Ouknine02b} D. Nualart and Y. Ouknine, Regularization of differential equations by fractional noise,  \textit{Stochastic Process. Appl.} {\bf 102} (2002), 103--116.


\bibitem{NR02} D. Nualart and A. R\u{a}\c{s}canu, Differential equations driven by fractional Brownian motion, \textit{Collect. Math.} {\bf 53} (2002), 55--81.

\bibitem{Nualart&Saussereau09} D. Nualart and B. Saussereau, Malliavin calculus for stochastic differential equations driven by a fractional Brownian motion, \textit{Stochastic Process. Appl.} {\bf 119} (2009), 391--409.




\bibitem{RW} P. Ren and F.-Y. Wang, Bismut formula for Lions derivative of distribution dependent SDEs and applications, \textit{J. Differential Equations} {\bf 267} (2019), 4745--4777.



\bibitem{SKM93} S. G. Samko, A. A. Kilbas and O. I. Marichev, Fractional Integrals and Derivatives, Theory and Applications, Gordon and Breach Science Publishers, Yvendon, 1993.

\bibitem{Song} Y. Song,  Gradient estimates and exponential ergodicity for Mean-Field SDEs with jumps,
 \textit{J. Theoret. Probab.} {\bf 33} (2020), 201--238.

\bibitem{Stein70} E. M. Stein, Singular Integrals and Differentiability Properties of Functions, Princeton University Press, Princeton, 1970.

\bibitem{SY21} Y. Suo and C. Yuan, Central limit theorem and moderate deviation principle for McKean-Vlasov SDEs,
\textit{Acta Appl. Math.} {\bf 175} (2021), 1--19.




\bibitem{Wang18} F.-Y. Wang,  Distribution dependent SDEs for Landau type equations, \textit{Stochastic Process. Appl.} {\bf 128} (2018), 595--621.

\bibitem{W22} F.-Y. Wang, Killed distribution dependent SDE for nonlinear Dirichlet problem, \textit{arXiv:2206.09115}.

\bibitem{WZZ15} R. Wang, J. Zhai and T. Zhang, A moderate deviation principle for 2-D stochastic Navier-Stokes equations, \textit{J. Differential Equations} {\bf 258} (2015), 3363--3390.




\bibitem{Zahle} M. Z\"{a}hle, Integration with respect to fractal functions and stochastic calculus I, \textit{Probab. Theory Related Fields} {\bf 111} (1998), 333--374.

\bibitem{ZZ15} J. Zhai and T. Zhang, Large deviations for 2-D stochastic Navier-Stokes equations driven by multiplicative l\'{e}vy noises, \textit{Bernoulli} {\bf 21} (2015), 2351--2392.

\bibitem{Zhang09} X. Zhang, A variational representation for random functionals on abstract Wiener spaces, \textit{J. Math. Kyoto Univ.} {\bf 49} (2009), 475--490.








}\end{thebibliography}
\end{document}